\documentclass[oneside]{amsart}
\usepackage{relsize}
\usepackage[bbgreekl]{mathbbol}
\usepackage{amsfonts}

\DeclareSymbolFontAlphabet{\mathbb}{AMSb} 
\DeclareSymbolFontAlphabet{\mathbbl}{bbold}
\newcommand{\Prism}{{\mathlarger{\mathbbl{\Delta}}}}
\usepackage{geometry}
\usepackage{needspace}
\RequirePackage{amsmath,amssymb,amsthm, amscd, comment,mathtools}
\usepackage{kantlipsum}
\usepackage[all]{xy}
\usepackage{graphicx}
\usepackage{tikz-cd}
\usepackage{rotating}
\usepackage{mathrsfs}
\usepackage[colorlinks=true, linkcolor=blue, urlcolor=blue, citecolor=blue]{hyperref}
\newtheorem{theorem}[subsection]{Theorem}
\newtheorem{proposition}[subsection]{Proposition}
\newtheorem{corollary}[subsection]{Corollary}
\newtheorem{lemma}[subsection]{Lemma}

\theoremstyle{definition}
\newtheorem{definition}[subsection]{Definition}
\newtheorem{remark}[subsection]{Remark}

\newtheorem{example}[subsection]{Example}
\newtheorem{assumption}[subsection]{Assumption}
\numberwithin{equation}{section}

\begin{document}
	\title{Deformations of Calabi--Yau varieties over transversal or crystalline prisms}
	\author{Tingxiao Shi}
	\begin{abstract}
	Recent work of Brantner and Taelman shows that derived deformation theory and period maps can recover strong mixed-characteristic lifting results. In this paper, we develop a prismatic version of this method for Calabi--Yau varieties over transversal or crystalline prisms, constructing compatible deformations over $A/I^n$ for every $n$ under suitable hypotheses.
	\par We also apply our method to $F$-split Calabi--Yau varieties over a field $k$ of characteristic $p>2$. In particular, we consider an $F$-split smooth proper variety $X$ of dimension $d$ over an algebraically closed field $k$ of characteristic $p>2$, with $H^{d-1}(X,O_X)=0$ and trivial canonical bundle. We show that $X$ admits a lift to $W_3(k)$.
	\end{abstract}

	\maketitle
	\setcounter{tocdepth}{1}
	\tableofcontents
	\section{Introduction}
	The Bogomolov--Tian--Todorov theorem asserts that smooth proper Calabi--Yau varieties over a field of characteristic zero have unobstructed formal deformations, but in positive characteristics, this statement no longer holds. Nevertheless, recent work shows that deformation-theoretic methods still gives lifting results under suitable hypotheses. In particular, Brantner and Taelman \cite{CY} prove mixed-characteristic unobstructedness and liftability results for Calabi--Yau varieties in characteristic $p$ by comparing the deformation functor of $X$ with a derived period functor. Some related theories also appear in the Serre--Tate theory of Calabi--Yau varieties developed by Achinger and Zdanowicz \cite{ST}, where the Frobenius splitting and the Hodge filtration gives the smooth deformation space under certain conditions.
	
	The purpose of this article is to develop a prismatic version of this deformation theory. Let $(A,I)$ be a transversal or crystalline prism. We study deformations of a smooth proper $p$-adic formal scheme $X$ over $A/I$ with trivial canonical bundle. The guiding principle is that the deformation problem of $X$ can be tested against a period-type deformation problem constructed from prismatic cohomology and the $p$-completed divided power de Rham complex $dR^{pd}_{X/(A/I)}$.
	
	An important motivating case is $(A,I)=(A_{\inf},(\xi))$, where $C=\mathbb C_p$, $A_{\inf}=W(\mathcal O_C^{\flat})$, and $(\xi)=\ker(\theta:A_{\inf}\to\mathcal O_C)$. A chosen first-order lift of a smooth formal scheme over $\mathcal O_C$ to $A_{\inf}/\xi^2$ is a standard geometric input in lift-dependent constructions of the $p$-adic Simpson correspondence. For example, see \cite[Theorem 1.1]{WangSimpson} for the correspondence between small generalised representations and small Higgs bundles. Compatible higher liftings also enter geometric $p$-adic Riemann--Hilbert theory: in \cite[Theorem 1.1]{LMNQW}, Liu--Ma--Nie--Qu--Wang construct correspondences between Hitchin-small $\mathbb B_{\mathrm{dR}}^+/t^n$-local systems and Hitchin-small integrable $t$-connections on the induced de Rham thickenings . Here $B_{\mathrm{dR}}^+$ is the $\xi$-adic completion of $A_{\inf}[1/p]$, so $B_{\mathrm{dR}}^+/\xi^n\simeq(A_{\inf}/\xi^n)[1/p]$, and $(t)=(\xi)$ in $B_{\mathrm{dR}}^+$. The integral lifting problem studied here thus provides geometric input for these correspondences. Our results give sufficient conditions for its solvability in the Calabi--Yau setting.
	\subsection*{Main results}
	This article has two main contributions. First, the method of \cite{CY}, which is formulated over Witt vectors in characteristic $p$, extends to general transversal or crystalline prisms.
	Second, we refine the period functor. Instead of working only with the full period functor $Per_X$ used in \cite{CY}, we introduce a smaller period functor $Per^1_X$ attached to the composite
	\[
	R\Gamma(dR^{pd}_{X/(A/I)})\longrightarrow R\Gamma(X,O_X)\longrightarrow H^d(X,O_X)[-d].
	\]
	Thus $Per^1_X$ remembers only the deformation of the top Hodge quotient relevant to Calabi--Yau obstruction theory. This modification is technically useful: it reduces the freeness hypotheses needed in the obstruction argument from all cohomological degrees to the range controlled by the Tor-dimension of $A$.
	
	Our main theorem is the following liftability result.
\par\medskip
	\noindent\textbf{Theorem A.} \textit{Let $A$ be an integral ring, and let $(A,I)$ be a transversal prism (or a crystalline prism with $p>2$). Let $m$ be the Tor-dimension of $A$, and let $X$ be a smooth proper $p$-adic formal scheme over $A/I$ of dimension $d$ with trivial canonical bundle. Assume:
	\\(a) The Hodge--de Rham spectral sequence of $X$ degenerates at $E_1$.
	\\(b) $H^i_{\Prism}(X)$ is free for all $i\in[d,d+m]$.
	\\Then there exists a compatible system of smooth proper $p$-adic formal deformations $X_n$ over $A/I^n$, for all $n\geqslant1$, with $X_1\cong X$. This system defines a smooth proper formal deformation $\mathscr{X}$ of $X$ over $\operatorname{Spf}_{(p,I)}A$.}
	\par\noindent The crystalline case at $p=2$ is discussed in Remarks \ref{char2}; in particular, the same conclusion holds for $(W(k),(2))$ over a perfect field $k$ without an additional lifting hypothesis.
\par\medskip
	The range $[d,d+m]$ is a consequence of the passage from $Per_X$ to $Per^1_X$. Under the above hypotheses, prismatic Poincar\'e duality and the Tor-amplitude of $A$ imply that the part of the dual prismatic complex which maps to the top Hodge quotient has the required freeness. This suffices to make the transition maps of the period functor $Per^1_X$ surjective and hence to eliminate obstructions.
	
	We also obtain an application to $F$-split Calabi--Yau varieties. Specializing to the crystalline prism $(W(k),p)$, the results of Achinger and Zdanowicz \cite{ST} identify Frobenius splitting with height one of the Artin--Mazur formal group and give a Frobenius-stable decomposition of crystalline cohomology. Together with Hodge--de Rham degeneration, which follows from Petrov's theorem \cite[Theorem 1.1]{Pet25}, this verifies the assumptions needed for the refined period map. We prove that if a good lift (in the sense that the first step of the Hodge filtration is preserved by Frobenius) of an $F$-split Calabi--Yau variety exists over $W_n(k)$, then it can be lifted further to $W_{n+1}(k)$. In particular, we obtain the following result.
\par
	\noindent\textbf{Theorem B.} Let $X$ be an $F$-split smooth proper variety of dimension $d$ over an algebraically closed field $k$ of characteristic $p>2$, with trivial canonical bundle and $H^{d-1}(X,O_X)=0$. Then $X$ admits a lift to $W_3(k)$.
	\subsection*{Sketch of proof}
	We briefly outline the proof of Theorem~A.
	We first construct the divided power deformation functor $Def_X^{pd}$ and the period map
	\[
		Def_X\longrightarrow Per_X.
	\]
 On shifted square-zero extensions, we compute the tangent map using the natural morphism
	\[
		Fil^1dR_{X/(A/I)}\longrightarrow L_{X/(A/I)}[-1].
	\]
	For the discrete thickenings needed in the obstruction argument, we then construct a canonical comparison map from deformations of $X$ to $Per^1_X$ and show that it is compatible with the period map above.
	
	After composing with the projection to $H^d(X,O_X)[-d]$, Serre duality identifies the obstruction map with the dual of the usual contraction pairing. Under the hypotheses of Theorem A, this gives the required injectivity on obstruction spaces. For a transversal prism, the first step uses $A/I^2$ with trivial higher divided powers on $I/I^2$: this is a PD elementary extension of $A/I$, so the period argument supplies the initial lift. For a crystalline prism with $p>2$, the canonical divided powers give the same first step. The deformation problem is then treated by induction along the thickenings $A/I^{n+1}\to A/I^n$, using the closed ordinary powers of the PD ideal for $n\geqslant2$. Injectivity of the period map on obstruction spaces forces each obstruction to vanish. Taking the limit of the structure sheaves gives a formal deformation over $\operatorname{Spf}_{(p,I)}A$.
	
	The article is organized as follows. Section~2 recalls divided power rings, animated divided power rings, and the $p$-completed divided power de Rham complex. Section~3 constructs the period functors $Per_X$ and $Per^1_X$. Section~4 proves the liftability theorem over transversal or crystalline prisms. Section~5 specializes to the crystalline prism $(W(k),p)$ and studies applications to $F$-split Calabi--Yau varieties and their canonical lifts.
	\subsection*{Conventions}
		Throughout this article, we assume that $A$ is an integral ring and that $(A,I)$ is a transversal or crystalline prism, as defined in \cite{prism} and \cite{absolute}. Let $(B,J,\gamma)$ be the $p$-adically complete PD envelope of $(A,I)$. To simplify notation, when $(A,I)=(A,p)$ is crystalline, we continue to use expressions such as ``formal schemes over $A/I$'' and ``$p$-adically completed modules'', which in this case simply mean schemes and modules over $A/I$, respectively. We frequently use the following assumption:
		\begin{assumption}\label{ass}
		$X$ is a smooth proper $p$-adic formal scheme over $A/I$ of relative dimension $d$, with trivial canonical bundle. Fix an identification $H^d(X,O_X)\cong A/I$, and suppose that $X$ satisfies the following conditions:
		\begin{center}
			(a) The canonical composite $q:\varphi^*R\Gamma_{\Prism}(X)\rightarrow R\Gamma(X,O_X)\rightarrow H^d(X,O_X)[-d]\cong(A/I)[-d]$ lifts through $A[-d]\to(A/I)[-d]$ to a map $\lambda:\varphi^*R\Gamma_{\Prism}(X)\to A[-d]$ admitting a section. We choose such a lift and section, giving a compatible splitting $\varphi^*R\Gamma_{\Prism}(X)\cong A[-d]\oplus M$.
			\\(b) $(\tau_{\geqslant 1-d}\varphi^*R\Gamma_{\Prism}(X))\hat{\otimes}^L_A A/I\rightarrow \tau_{\geqslant1-d}R\Gamma(X,O_X)$ is surjective in every degree.
		\end{center}
		Here, $\varphi$ is the Frobenius action on the $\delta$-ring $A$, and $\varphi^*R\Gamma_{\Prism}(X):=R\Gamma_{\Prism}(X)\hat{\otimes}^L_{A,\varphi} A$.
	\end{assumption}
	\subsection*{Acknowledgements}
		The author thanks Yupeng Wang for suggesting this research topic and giving advice on the research direction; and also thanks my tutor Ruochuan Liu for helpful advice on the structure of the article, and useful suggestion on giving examples of the paper.
	\section{Divided power rings and de Rham complexes}
	To prove our main theorem, we introduce the derived divided power de Rham complex, following the discussion in \cite{mao}. We begin with the definition of a divided power ring.
	\begin{definition}(Divided power rings, \cite[\href{https://stacks.math.columbia.edu/tag/07GL}{07GL}]{stack})\label{pdring}
		Let $A$ be a ring, and let $I$ be an ideal of $A$. A collection of maps $\gamma_n:I\rightarrow I, n>0$ is called a \textit{divided power structure} on $I$ if, for all $n\geqslant0,\  m>0,\ x,y\in I,\ a\in A$, the following identities hold:
		\\(1) $\gamma_0(x)=1,\ \gamma_1(x)=x$;
		\\(2) $\gamma_n(x)\gamma_m(x)=\frac{(m+n)!}{m!n!}\gamma_{m+n}(x)$;
		\\(3) $\gamma_n(ax)=a^n\gamma_n(x)$;
		\\(4) $\gamma_n(x+y)=\sum\limits_{i=0}^n\gamma_i(x)\gamma_{n-i}(y)$;
		\\(5) $\gamma_n(\gamma_m(x))=\frac{(nm)!}{n!(m!)^n}\gamma_{mn}(x)$.
		\\A \textit{divided power ring} is a triple $A=(A,I,\gamma)$. A morphism of divided power rings is a morphism of rings $f:A\rightarrow A'$ mapping $I\rightarrow I'$ and satisfying $\gamma_n'(fx)=f\gamma_n(x)$. We denote the category of divided power rings by $PDRing$.
		
	\end{definition}
	\begin{definition}(Animated divided power rings)
		For a ring $R$ and an $R$-module $M$, we define the \textit{free divided power algebra} $\Gamma_R(M)$ as the quotient of the free commutative $R$-algebra on symbols
		\[
		\{\gamma_n(y)|y\in M,\ n\geqslant 0\}
		\]
		by the standard divided power relations. Equivalently, $\Gamma_R(M)$ is characterized by the natural bijection between morphisms of divided power $R$-algebras $\Gamma_R(M)\to(S,K,\delta)$ and $R$-linear maps $M\to K$. Write $PDPoly\subset PDRing$ for the full subcategory spanned by all divided power rings of the form
		\[
		\Gamma_{\mathbb{Z}[x_1,\cdots,x_n]}(y_1,\cdots,y_m),\qquad n,m\geqslant0.
		\]
		Write $ArrPDPoly$ for the category whose objects are the standard inclusions obtained by adjoining finitely many polynomial and divided power generators,
		\[
		\Gamma_{\mathbb{Z}[x_1,\ldots,x_n]}(y_1,\ldots,y_m)
		\longrightarrow
		\Gamma_{\mathbb{Z}[x_1,\ldots,x_{n+n'}]}(y_1,\ldots,y_{m+m'}),
		\]
		and whose morphisms are commutative squares of divided power rings; see \cite[Definition 4.22 and Lemma 4.23]{CY}.
		\\We define the $\infty$-category of \textit{animated divided power rings} $CAlg^{pd}:=\mathscr{P}_{\Sigma}(PDPoly)$ by freely adjoining sifted colimits to $PDPoly$.
	\end{definition}
	Every ring $R$ has a trivial divided power structure $(R,0,0)$, which induces a natural map $CAlg^{an}\hookrightarrow CAlg^{pd}$. In the other direction, we have two maps from $CAlg^{pd}$ to $CAlg^{an}$:
	\begin{definition}
		Consider the following two functors from $PDPoly$ to $CAlg$:
		\[
		und:(A,I,\gamma)\mapsto A
		\]
		\[
			\overline{(\cdot)}:(A,I,\gamma)\mapsto A/I
		\]
		By taking sifted colimits, these two maps extend to maps $CAlg^{pd}\rightarrow CAlg^{an}$, which we continue to denote by $und$ and $\overline{(\cdot)}$, respectively.
	\end{definition}
	In \cite{Bhatt}, the author introduced the derived de Rham complex for morphisms of animated rings and, in Chapter 7, explained its $p$-adic completion; in \cite{mao}, the author extended the definition to the divided power de Rham complex over animated divided power rings.
	
	In this article, we combine these definitions, as in \cite[Chapter 4]{CY}, to define the $p$-completed divided power de Rham complex.
	\begin{definition}(\cite[\href{https://stacks.math.columbia.edu/tag/07HQ}{07HQ}]{stack})\label{uncomdR}
		For a morphism of divided power rings $(A,I,\gamma)\rightarrow (B,J,\delta)$ and a $B$-module $M$, a \textit{divided power $A$-derivation into $M$} is an additive map $\theta:B\rightarrow M$ that annihilates the elements of $A$, satisfies the Leibniz rule $\theta(bb')=b\theta(b')+b'\theta(b)$, and satisfies
		\[
		\theta(\delta_n(x))=\delta_{n-1}(x)\theta(x)
		\]
		for $n\geqslant 1$ and $x\in J$.
		\\In this setting, as for ordinary derivations, there exists a \textit{universal divided power $A$-derivation}
		\[
		d_{B/A,\delta}:B\rightarrow\Omega_{B/A,\delta}
		\]
		where $\Omega_{B/A,\delta}$ is the \textit{module of PD-differentials}: every divided power $A$-derivation $B\to M$ factors uniquely through $d_{B/A,\delta}$ by a $B$-linear map $\Omega_{B/A,\delta}\to M$.
		\\We can use this to construct the (uncompleted) divided power de Rham complex over free divided power algebras
		\[
		(\Omega^{\cdot}_{B/A,\delta},d)
		\]
		with terms $\Omega^i_{B/A,\delta}=\wedge^i\Omega^1_{B/A,\delta}$ and differential $d:\Omega^i_{B/A,\delta}\rightarrow\Omega^{i+1}_{B/A,\delta}$ defined by
		\[
		d(f_0df_1\wedge\cdots\wedge df_i)=df_0\wedge\cdots\wedge df_i
		\]
	\end{definition}
	\begin{definition}($p$-completed divided power de Rham complex)
		The construction in Definition \ref{uncomdR} gives a map
		\[
		\Omega^{\bullet}_{B/A,\delta}:ArrPDPoly\rightarrow CAlg^{pd}Mod:=CAlg^{pd}\times_{CAlg^{an}} CAlg^{an}Mod
		\]
		We define the $\infty$-category $CAlg^{pd}Mod^\wedge$ as the full subcategory of $CAlg^{pd}Mod$ spanned by all derived $p$-complete module complexes.
		\\Then we have a natural completion map
		\[
		CAlg^{pd}Mod\rightarrow CAlg^{pd}Mod^\wedge
		\]
		which is left adjoint to the inclusion map.
		\\Combining these maps and taking sifted colimits, we construct the \textit{$p$-completed divided power de Rham complex}
		\[
		dR^{pd}_{\bullet/\bullet}:Fun(\Delta^1,CAlg^{pd})\rightarrow  CAlg^{pd}Mod^{\wedge}
		\]
		as the unique sifted-colimit-preserving extension of the composition of $\Omega^{\bullet}_{B/A,\delta}$ and the $p$-adic completion map.
		Similarly, we write $L^{pd}_{B/A}$ for the sifted-colimit-preserving extension of the functor $\Omega^1_{B/A,\delta}$.
	\end{definition}
	\begin{remark}
		(1) In this article, we follow the notation of \cite{CY}: our $dR^{pd}_{B/A}$ is the $p$-adic completion of the de Rham complex defined in \cite{mao}. In what follows, we still invoke some of the results in \cite{mao}, but always after passing to their $p$-completed versions.
		\\(2) As noted in \cite[Remark 4.5]{mao}, although the module of PD-differentials is independent of the PD-structure on $A$, the
		derived divided power de Rham complex does depend on the PD-structure on $A$.
	\end{remark}
	This definition extends readily to the setting of sheaves. For any site $\mathscr{T}$ with a sheaf $O_X\in Shv(\mathscr{T},CAlg^{pd})$ over an animated divided power ring $A$, we write $dR^{pd}_{O_X/A}$ for the sheafification of the presheaf $\{U,dR^{pd}_{O_X(U)/A}\}$.
	
	Like the usual derived de Rham complex, the divided power de Rham complex has the following basic properties:
	\begin{proposition}\label{basicprop}
		(1) (Base change) For a pushout $B'=B\otimes_A A'$ in $CAlg^{pd}$, the natural map $dR^{pd}_{B/A}\hat{\otimes}^L_A A'\rightarrow dR^{pd}_{B'/A'}$ is an equivalence.
		\\(2) (Filtration) The naive filtration
		\[
		Fil^n\Omega^\bullet_{B/A,\delta}=[\cdots\rightarrow 0\rightarrow \Omega^n_{B/A,\delta}\rightarrow \Omega^{n+1}_{B/A,\delta}\rightarrow\cdots]
		\]
		on $\Omega^\bullet_{B/A,\delta}$ induces a natural filtration on $dR^{pd}_{B/A}$, which satisfies
		\[
		gr^0dR^{pd}_{B/A}\simeq B^{\wedge}_p,\qquad gr^1dR^{pd}_{B/A}\simeq (L^{pd}_{B/A})^{\wedge}_p[-1].
		\]
		\\(3) (Comparison with de Rham cohomology) For any map $A\rightarrow B$ between animated rings (viewed as animated divided power rings with trivial divided powers), $dR^{pd}_{B/A}$ agrees with Bhatt's $p$-completed variant of Illusie's derived de Rham cohomology $(dR_{B/A})^{\wedge}_p$.
		\\(4) (Comparison with crystalline cohomology) Assume that $(B,J,\gamma)$ is a discrete PD-ring and $R$ is a smooth $B/J$-algebra. Then $dR^{pd}_{(R,0,0)/(B,J,\gamma)}\cong (R\Gamma_{cris}(R/B))^{\wedge}_p$.
	\end{proposition}
	\begin{proof}
		For (1)\&(2), these properties are easily verified for maps in $ArrPDPoly$; the sifted-colimit-preserving property then gives the desired results for maps in $Fun(\Delta^1,CAlg^{pd})$.
		\\For (3), see \cite[Lemma 4.10]{mao}.
		\\For (4), smoothness implies that $R$ is quasisyntomic over $B/J$. By \cite[Propositions 4.87 and 4.90]{mao}, the uncompleted divided power de Rham complex agrees with the cohomology of the affine crystalline site and then with classical crystalline cohomology. Taking derived $p$-completion gives the asserted comparison.
	\end{proof}
	We use the $p$-completed divided power de Rham complex because it exhibits a form of rigidity under deformations, much like the usual de Rham complex over the complex numbers.
	\begin{proposition}(\cite[Proposition 4.16]{mao})
		For every map $C\rightarrow D$ in $CAlg^{pd}$, the induced map
		\[
		dR^{pd}_{D/C}\rightarrow dR^{pd}_{\overline{D}/C}
		\]
		is an equivalence in $Mod^{\wedge}_C$.
	\end{proposition}
	\begin{proof}
		The uncompleted version of this proposition was proved in \cite[Proposition 4.16]{mao}. Taking $p$-adic completion gives the desired equivalence.
	\end{proof}
	This rigidity property helps us construct a canonical map from the sheaf of divided power de Rham complexes to the structure sheaf of deformations of the original formal scheme, but this is not working for any divided power rings.
	
	We first specify the animated divided power rings under consideration.
	\begin{definition}
		(1) For any $(R,I_R,\eta)\in PDPoly$ and any $R$-module $V$, we write $sqz^{pd}_{R}(V)=(R\oplus V,I_R\oplus V,\eta')$ for the square-zero extension of $R$ by $V$. Here the divided power structure $\eta'$ is defined by
		\[
		\eta'_n(x+y)=\eta_n(x)+\eta_{n-1}(x)y,\quad x\in I_R,\ y\in V,
		\] for $n\geqslant1$, with $\eta'_0(x+y)=1$.
		\\Then, by taking sifted colimits over square-zero extensions of modules over divided power polynomial rings, we obtain a square-zero extension map
		\[
		sqz^{pd}_{C}:Mod^{cn}_{C}\rightarrow CAlg^{pd}
		\]
		Here, $Mod^{cn}_{C}$ is the $\infty$-category of all connective $C$-modules.
		\\(2) We write $sqz_C$ for $und(sqz^{pd}_C): Mod_C^{cn}\rightarrow CAlg^{an}$. By \cite[Corollary 4.7]{CY}, $und$ preserves limits and small colimits, so this agrees with the square-zero extension for animated algebras.
	\end{definition}
	\begin{definition}\label{category}
		(1) We define $CAlg^{an}_{B//(A/I)}$ as the $\infty$-category consisting of all animated rings $N$ with a diagram
		\[
		B\rightarrow N\rightarrow A/I
		\]
		Similarly, $CAlg^{pd}_{B//(A/I)}$ is defined as the $\infty$-category consisting of all animated divided power rings $N$ with a diagram
		\[
		(B,J,\gamma)\rightarrow N\rightarrow (A/I,0,0)
		\]
		Here $(B,J,\gamma)$ denotes the $p$-complete PD-envelope of $(A,I)$.
		\\(2) We say that a map $P\rightarrow Q$ between two objects of $CAlg^{an}_{B//(A/I)}$ is an \textit{elementary extension} if there exists a discrete $A/I$-module $V$ such that, for some $n\geqslant 0$, the map fits into the following pullback square:
		\[
		\begin{tikzcd}
			P \arrow[d] \arrow[r] \arrow[dr, phantom, "\lrcorner", at start]& Q \arrow[d]  \\
			A/I \arrow[r]         & sqz_{A/I}(V[n])
		\end{tikzcd}
		\]
		\\Similarly, a map $P\rightarrow Q$ between two objects of $CAlg^{pd}_{B//(A/I)}$ is an \textit{elementary extension} if there exists a discrete $A/I$-module $V$ such that, for some $n\geqslant 0$, the map fits into the following pullback square:
		\[
		\begin{tikzcd}
			P \arrow[d] \arrow[r] \arrow[dr, phantom, "\lrcorner", at start]& Q \arrow[d]  \\
			A/I \arrow[r]         & sqz^{pd}_{A/I}(V[n])
		\end{tikzcd}
		\]
		(3) Let $CAlg^{pd,ele}_{B//(A/I)}\subset CAlg^{pd}_{B//(A/I)}$ be the full subcategory of diagrams $(B,J,\gamma)\rightarrow N\rightarrow(A/I,0,0)$ satisfying the following conditions:
		\\(1) $\pi_0(B)\rightarrow\pi_0(N)$ is surjective and $\overline{N} = A/I$.
		\\(2) The map $und(N)\rightarrow A/I$ can be decomposed into a finite sequence of elementary extensions.
		\\(3) The underlying animated ring $N$ is derived $p$-complete.
	\end{definition}
	\begin{proposition}\label{pullback}
		(1) For any discrete $A/I$-module $V$, there is a pullback square
		\[
		\begin{tikzcd}
			sqz^{pd}_{A/I}(V[n]) \arrow[d] \arrow[r] & A/I \arrow[d]     \\
			A/I \arrow[r]                    & {sqz^{pd}_{A/I}(V[n+1])}
		\end{tikzcd}
		\]
		for any $n\geqslant 0$.
		\par\Needspace{9\baselineskip}\noindent(2) Let $(B,J,\gamma)$ be the $p$-complete PD-envelope of $(A,I)$, with $F_n=\overline{J^n}^{\,p}$. For any $n>0$, there is a pullback square
		\[
		\begin{tikzcd}
			B/F_{n+1} \arrow[d] \arrow[r] & B/F_n \arrow[d]       \\
			A/I \arrow[r]                 & {sqz_{A/I}(F_n/F_{n+1}[1])}
		\end{tikzcd}
		\]
		Moreover, suppose that one of the following conditions holds:
		\\(a) $n\geqslant 2$;
		\\(b) $p>2$ and $(A,I)$ is crystalline.
		\\Then there is a pullback square
		\[
		\begin{tikzcd}
			(B/F_{n+1}, J/F_{n+1}, \gamma) \arrow[d] \arrow[r] & (B/F_n, J/F_n, \gamma) \arrow[d]       \\
			A/I \arrow[r]                 & {sqz^{pd}_{A/I}(F_n/F_{n+1}[1])}
		\end{tikzcd}
		\]
	\end{proposition}
	\begin{proof}
		(1) We construct a commutative diagram
		\[
		\begin{tikzcd}
			sqz^{pd}_{A/I}(V[n]) \arrow[d] \arrow[r] & A/I \arrow[d]     \\
			A/I \arrow[r]                    & {sqz^{pd}_{A/I}(V[n+1])}
		\end{tikzcd}
		\]
		induced by the bicartesian square
		\[
		\begin{tikzcd}
			V[n] \arrow[d] \arrow[r] & 0 \arrow[d] \\
			0 \arrow[r]           & {V[n+1]}     
		\end{tikzcd}
		\]
		in $Mod_A$.
		\\Therefore, by checking the fibers, we see that the first diagram is also a pullback square for every $n\geqslant0$. If $V$ is derived $p$-complete and $n\geqslant1$, this also gives $sqz^{pd}_{A/I}(V[n])\in CAlg^{pd,ele}_{B//(A/I)}$. For $n=0$ and $V\ne0$, the structural map $B\rightarrow A/I\oplus V$ is not surjective on $\pi_0$, so this membership assertion does not apply.
		\\(2) In the transversal case, $B/J=A/I$ is $p$-torsion-free, so $J\cap p^rB=p^rJ$. If $x-y=p^rz$ with $x,y,z\in J$, the divided power addition formula gives
		\[
		\gamma_m(x)-\gamma_m(y)=\sum_{k=1}^m p^{rk}\gamma_{m-k}(y)\gamma_k(z)\in p^rB.
		\]
		Thus each $\gamma_m$ is continuous on $J$. The ordinary powers $J^n$ are PD-stable by the divided power product and addition formulas, so their closures $F_n$ are PD-stable as well. In the crystalline case, $B=A$, $J=(p)$, and $F_n=(p^n)$ with the canonical divided powers, so PD-stability also holds.
		
		Now \cite[Proposition
		25.3.6.1]{SAG} gives an isomorphism
		\[
		F_n/F_{n+1}\cong \pi_1 L_{(B/F_n)/(B/F_{n+1})}
		\]
		Here, $L_{(B/F_n)/(B/F_{n+1})}$ is the cotangent complex.
		\\This induces a map
		\[
		L_{(B/F_n)/(B/F_{n+1})}\rightarrow \tau_{\leqslant 1}L_{(B/F_n)/(B/F_{n+1})}\rightarrow F_n/F_{n+1}[1]
		\]
		which therefore gives a map
		\[
		B/F_n\rightarrow sqz_{B/F_n}(F_n/F_{n+1}[1])\rightarrow sqz_{A/I}(F_n/F_{n+1}[1])
		\]
		Hence, we have a commutative diagram
		\[
		\begin{tikzcd}
			B/F_{n+1} \arrow[d] \arrow[r] & B/F_n \arrow[d]       \\
			A/I \arrow[r]                 & {sqz_{A/I}(F_n/F_{n+1}[1])}
		\end{tikzcd}
		\]
		Checking the fibers shows that this is a pullback square.
		\\Similarly, for PD-algebras, \cite[Proposition 4.39]{CY} gives a map
		\[
		L^{pd}_{(B/F_n)/(B/F_{n+1})}\rightarrow \tau_{\leqslant 1}L^{pd}_{(B/F_n)/(B/F_{n+1})}\cong F_n/(F_{n+1}+(F_n)^{[2]})[1]
		\]
		Here, $(F_n)^{[2]}$ denotes the ideal generated by $\gamma_m(x)$ with $m\geqslant2$ and $x\in F_n$.
		\\If $n\geqslant 2$, for any $ay\in J^n$ with $a\in J^{n-1}$ and $y\in J$, we have
		\[
		\gamma_m(ay)=a^m\gamma_m(y)\in J^{2n-1}\subset J^{n+1}
		\]
		The divided power addition formula extends this calculation to finite sums in $J^n$, since every mixed term contains at least two factors from $J^n$. In the transversal case, continuity then gives $\gamma_m(F_n)\subset F_{n+1}$. In the crystalline case, $F_n=J^n$, so the same inclusion already follows. Thus $(F_n)^{[2]}\subset F_{n+1}$ for $n\geqslant2$.
		\\On the other hand, if $I=(p)$ and $p>2$, then $B=A$ and $F_n=(p^n)$ with the canonical divided powers. The only additional case is $n=1$: for $m\geqslant2$, the inequality $m-v_p(m!)\geqslant2$ gives $\gamma_m(p)=p^m/m!\in(p^2)$. Thus $(F_n)^{[2]}\subset F_{n+1}$ also holds in this case.
		\\Thus, in both cases (a) and (b), we obtain a map
		\[
		B/F_n\rightarrow sqz^{pd}_{B/F_n}(F_n/F_{n+1}[1])\rightarrow sqz^{pd}_{A/I}(F_n/F_{n+1}[1])
		\]
		which gives a commutative diagram
		\[
		\begin{tikzcd}
			B/F_{n+1} \arrow[d] \arrow[r] & B/F_n \arrow[d]       \\
			A/I \arrow[r]                 & {sqz^{pd}_{A/I}(F_n/F_{n+1}[1])}
		\end{tikzcd}
		\]
		Checking the fibers of the horizontal maps in this diagram shows that it is a pullback square. Finally, the ordinary elementary tower $B/F_n\to\cdots\to B/F_1=A/I$, together with the completeness just proved, shows that $B/F_n\in CAlg^{pd,ele}_{B//(A/I)}$ for every $n\geqslant1$.
	\end{proof}
	\Needspace{9\baselineskip}
	\begin{lemma}[The first transversal thickening]\label{firstthickening}
		Let $(A,I)$ be a transversal prism. Equip $I/I^2\subset A/I^2$ with the divided powers $\eta_1(x)=x$ and $\eta_r(x)=0$ for $r\geqslant2$. Then $A/I^2\to A/I$ is a PD elementary extension. More precisely, there is a pullback square
		\[
		\begin{tikzcd}
			(A/I^2,I/I^2,\eta) \arrow[d] \arrow[r] \arrow[dr, phantom, "\lrcorner", at start] & A/I \arrow[d] \\
			A/I \arrow[r] & {sqz^{pd}_{A/I}(V[1])}.
		\end{tikzcd}
		\]
		Moreover, if $F'_2=\overline{J^{[2]}}^{\,p}$, then $B/F'_2\cong A/I^2$ as PD $A$-algebras.
	\end{lemma}
	\begin{proof}
		First, by
		\[
		0\longrightarrow I/I^2\longrightarrow A/I^2\longrightarrow A/I\longrightarrow0
		\]
		we can see that $A/I^2$ is $p$-torsion-free and $p$-complete. The stated divided powers satisfy the PD identities because $(I/I^2)^2=0$. They are also compatible with the canonical divided powers on $p$.
		\\The augmentation has kernel $I/I^2$, with $(I/I^2)^{[2]}=0$, thus \cite[Proposition 4.39]{CY} gives $\tau_{\leqslant1}L^{pd}_{R/C}\simeq V[1]$, which gives the pullback diagram above and hence prove that $(A/I^2,I/I^2,\eta)$ is in $CAlg^{pd,ele}_{B//R}$.
		\\Finally, the kernel of the map $B\rightarrow A/I^2$ contains $F'_2$; conversely,
		\[
		xy=\gamma_2(x+y)-\gamma_2(x)-\gamma_2(y)
		\]
		gives $J^2\subset J^{[2]}$, so $A\rightarrow B/F'_2$ factors through $A/I^2$. Therefore, the map $B/F'_2\rightarrow A/I^2$ is an isomorphism.
	\end{proof}
	In what follows, we mainly consider animated divided power rings in $CAlg^{pd,ele}_{B//(A/I)}$. In this category, we have the following proposition:
	\begin{proposition}\label{rigidity}
		Assume that $X$ is a smooth, proper $p$-adic formal scheme over $A/I$. For any $C\in CAlg^{pd,ele}_{B//(A/I)}$ and any $O^{pd}_\mathscr{X}$ satisfying $O^{pd}_\mathscr{X}\hat{\otimes}^L_{C} A/I=O_X$, we have a canonical isomorphism
		\[
		R\Gamma(X,dR^{pd}_{O^{pd}_{\mathscr{X}}/C})\cong R\Gamma_{cris}(X/B)\hat{\otimes}^L_B C
		\]
		Here, $(B,J,\gamma)$ is the $p$-complete PD-envelope of $(A,I)$. The assertion also applies to split square-zero objects $C=sqz^{pd}_{A/I}(V[n])$, where $n\geqslant0$ and $V$ is discrete and derived $p$-complete.
	\end{proposition}
	\begin{proof}
		By Proposition \ref{basicprop}, we have
		\[
		\begin{split}
		R\Gamma_{cris}(X/B)\hat{\otimes}^L_B C&\cong R\Gamma(X,dR^{pd}_{O_X/(B,J,\gamma)})\hat{\otimes}^L_B C
		\\&\cong R\Gamma(X,dR^{pd}_{O_X\otimes^L_B C/C})
		\\&\cong R\Gamma(X,dR^{pd}_{\overline{O_X\otimes^L_B C}/C})
		\\&\cong R\Gamma(X,dR^{pd}_{O_X/C})
		\\&\cong R\Gamma(X,dR^{pd}_{O^{pd}_{\mathscr{X}}/C})
		\end{split}
		\]
		Here we also use the fact that the map $\overline{(\bullet)}$ preserves colimits. The proof uses the identity $\overline C=A/I$ and the base-change and PD-reduction properties, but not surjectivity of $B\to C$ on $\pi_0$; hence it also applies to the stated split square-zero objects.
	\end{proof}
	\begin{proposition}\label{conservative}
		For $P\in CAlg^{pd,ele}_{B//(A/I)}$ and any $M\in\widehat D(P)$, if $M\hat{\otimes}_P^L A/I$ is connective, then $M$ is connective.
	\end{proposition}
	\begin{proof}
		First, since $P\rightarrow A/I$ can be decomposed into a finite sequence of elementary extensions, $P$ has only finitely many nonzero homotopy groups, and $\pi_0(P)\rightarrow R$ has a nilpotent kernel $K$.
		
		The finite Postnikov filtration of $P$, as an ordinary $P$-module, has factors $\pi_i(P)[i]$ with $i\geqslant0$. Filtering each $\pi_i(P)$ by powers of $K$ refines it to a finite filtration whose factors are $W[i]$ for discrete $A/I$-modules $W$. For each such factor,
		\[
		M\hat{\otimes}_P^L W[i]\simeq (M\hat{\otimes}_P^L A/I)\hat{\otimes}_{A/I}^L W[i]
		\]
		Since derived $p$-completion preserves connective complexes, $M\hat{\otimes}_P^L W[i]$ is connective. The completed derived tensor product is exact, so the finite filtration shows that $M\hat{\otimes}_P^L P\simeq M$ is connective.
	\end{proof}
	\section{Period maps}
		Having recalled the definition and basic properties of animated divided power rings, we now construct the period map. We first define the deformation functors used in this construction.
	\begin{definition}\label{deformation}
			(1) (PD module deformations) Let $\mathscr{C}$ be an $A$-linear stable $\infty$-category. Consider the $\infty$-category
		\[
			\mathscr{D}=CAlg_A^{pd}\times _{CAlg^{an}_A} Mod(\mathscr{C}).
		\]
			Let $q:\mathscr{D}\rightarrow CAlg_A^{pd}$ be the projection, and let $\mathscr{D}^{coca}$ be the subcategory of $\mathscr{D}$ containing all objects and only the $q$-cocartesian morphisms.
		\\For any $A/I$-module $M\in Mod_{A/I}(\mathscr{C})$, we have a left fibration $\mathscr{D}^{coca}_{/(A/I,M)}\rightarrow CAlg_{A//(A/I)}^{pd}$, where the base records both the $A$-algebra structure and the augmentation to $A/I$.
			\\Straightening yields a functor $CAlg_{A//(A/I)}^{pd}\rightarrow\hat{S}$, called the deformation functor of $(A/I,M)$ and denoted by $Def_M^{pd}$.
			\\(2) (PD algebra deformations) For a topological space $\mathscr{T}$, consider the $\infty$-category
		\[
			\mathscr{D}=CAlg_A^{pd}\times_{Shv(\mathscr{T},CAlg^{pd})^{\wedge}}Fun(\Delta^1,Shv(\mathscr{T},CAlg^{pd})^{\wedge}).
		\]
			Here the map $CAlg_A^{pd}\rightarrow Shv(\mathscr{T},CAlg^{pd})^{\wedge}$ is obtained by taking the constant sheaf followed by derived $p$-adic completion.
			Similarly, we obtain a left fibration $\mathscr{D}^{coca}_{/(A/I,\mathscr{F}^{pd})}\rightarrow CAlg_{A//(A/I)}^{pd}$ for any sheaf $\mathscr{F}^{pd}$ over $A/I$.
			\\This defines the deformation functor $Def^{pd}_{\mathscr{F}}$.
	\end{definition}
		Now let $X=(\mathscr{T},O_X)$ be a smooth proper $p$-adic formal scheme over $A/I$. We define $O_X^{pd}\in Shv(\mathscr{T},CAlg^{pd})^{\wedge}$ using the natural injection $CAlg^{an}\rightarrow CAlg^{pd}$. Here $Shv(\mathscr{T},CAlg^{pd})^{\wedge}$ denotes the category of derived $p$-complete sheaves on $\mathscr{T}$.
	
	We will use the following deformation functors:
		\\$Def^{pd}_X$: the deformation functor of $(A/I,O^{pd}_X)$ in $Shv(\mathscr{T},CAlg^{pd})^{\wedge}$;
		\\$Def^{pd}_{R\Gamma(\theta_X)}$: the deformation functor of $R\Gamma(dR^{pd}_{X/(A/I)})\rightarrow R\Gamma(X,O^{pd}_X)$ in $Fun(\Delta^1, CAlg^{pd}Mod^{\wedge})$;
		\\$Def^{pd}_{R\Gamma(dR^{pd}_{X/(A/I)})} $: the deformation functor of $(A/I,R\Gamma(dR^{pd}_{X/(A/I)}))$ in $CAlg^{pd}Mod^{\wedge}$.
		\begin{definition}(Deformations of $O_X$)
			(1) For a smooth proper $p$-adic formal scheme $X=(\mathscr{T},O_X)$ over $A/I$, consider the $\infty$-category
		\[
		\mathscr{D}=CAlg^{an}_A\times_{Shv(\mathscr{T},CAlg^{an})^{\wedge}}Fun(\Delta^1,Shv(\mathscr{T},CAlg^{an})^{\wedge})
		\]
			without any PD structure. As above, this construction gives a left fibration
		\[
			\mathscr{D}^{coca}_{/(A/I,O_X)}\rightarrow CAlg^{an}_{A//(A/I)},
		\]
		which defines a functor $CAlg^{an}_{A//(A/I)}\rightarrow \hat{S}$, denoted by $Def_X$.
		\\By precomposing with the forgetful functor $und$, we also obtain $Def_X: CAlg^{pd}_{A//(A/I)}\rightarrow \hat{S}$. We use the same notation for these two functors.
			\\(2) For any $B\in CAlg^{pd}$, a sheaf $\mathscr{F}\in Shv(\mathscr{T},CAlg^{an})^{\wedge}$ over $und(B)$ is automatically a sheaf in $Shv(\mathscr{T},CAlg^{an})^{\wedge}$ with trivial PD structure. Thus the base-change functor $\bullet\hat{\otimes}^L_{und(B)}B$ in $Shv(\mathscr{T},CAlg^{pd})^{\wedge}$ sends a sheaf $\mathscr{F}$ over $und(B)$ to a sheaf $\mathscr{F}\hat{\otimes}^L_{und(B)}B$ over $B$.
		\\Since $und$ preserves small colimits by \cite[Corollary 4.7]{CY}, this is equivalent to giving $\mathscr{F}$ a PD structure.
			\\This construction therefore induces a functor $Def_X\rightarrow Def^{pd}_X$.
	\end{definition}
	\begin{remark}\label{classicalthickenings}
			(1) Unraveling the definition, we see that, set-theoretically, for any $C\in CAlg^{pd}$, an element of $Def_X^{pd}(C)$ corresponds to a sheaf $\mathscr{F}\in Shv(\mathscr{T},CAlg^{pd})^{\wedge}$ together with a cocartesian diagram:
		\[
		\begin{tikzcd}
			C \arrow[r] \arrow[d] & \mathscr{F} \arrow[d] \\
			A/I \arrow[r]         & O_X^{pd} \arrow[ul, phantom, "\ulcorner", very near start]            
		\end{tikzcd}\]
			The other deformation functors admit similar informal descriptions.
			\par (2) Assume $X$ is smooth and proper over $R=A/I$. For the discrete bases $C=A/I^n$ used in the lifting argument, a sheaf $\mathscr{F}\in Def_X(C)$ determines a smooth proper $p$-adic formal scheme $(|X|,\pi_0\mathscr{F})$. We justify this at finite level. In the transversal case, $C$ and $R$ are $p$-torsion-free. Set $C_r=C/p^r$, $R_r=R/p^r$, $X_r=X\times_R R_r$, and $\mathscr{F}_r=\mathscr{F}\otimes_C^L C_r$. Completed derived base change gives
		\[
		\mathscr{F}_r\otimes_{C_r}^L R_r\simeq O_{X_r}.
		\]
			The kernel of $C_r\to R_r$ is nilpotent. Since $O_{X_r}$ is flat over $R_r$, the nilpotent criterion for Tor-amplitude \cite[\href{https://stacks.math.columbia.edu/tag/0H75}{0H75}]{stack}, applied stalkwise, shows that $\mathscr{F}_r$ is discrete and flat over $C_r$. The successive quotients of its nilpotent ideal filtration are quasi-coherent, being obtained by tensoring the corresponding base-ideal quotients with $O_{X_r}$. Applying \cite[\href{https://stacks.math.columbia.edu/tag/05YV}{05YV}]{stack} iteratively therefore realizes $(|X|,\mathscr{F}_r)$ as a scheme; smoothness and properness follow from \cite[\href{https://stacks.math.columbia.edu/tag/06AG}{06AG}]{stack}.
			\par On an affine open $U$ of $X/p$, these thickenings are affine by \cite[\href{https://stacks.math.columbia.edu/tag/06AD}{06AD}]{stack}. Their rings $B_r=\Gamma(U,\mathscr{F}_r)$ satisfy $B_{r+1}/p^rB_{r+1}=B_r$, so $B_U=\varprojlim_r B_r$ is $p$-complete with $B_U/p^rB_U=B_r$ \cite[\href{https://stacks.math.columbia.edu/tag/09B8}{09B8}]{stack}. The surjective transition maps give $R\Gamma(U,\mathscr{F})\simeq B_U$. Sections on distinguished opens are the corresponding completed localizations, identifying $\mathscr{F}|_U$ with $O_{\operatorname{Spf}_p B_U}$. These charts glue via the restrictions of $\mathscr{F}$. The glued reductions $(|X|,\mathscr{F}_r)$ are smooth and proper, yielding the claimed smooth proper formal scheme. In the crystalline case, $p$ is nilpotent in $C=A/p^n$, so the same nilpotent argument applies directly over $C\to R$, without taking reductions modulo $p^r$.
	\end{remark}
		We will use the fact that these deformation functors preserve pullback squares under the following $\pi_0$-surjectivity hypotheses.
	\begin{proposition}\label{formalmoduli}
			Suppose that the category $\mathscr{C}$ admits a left-complete $t$-structure and that $M\in Mod_{A/I}(\mathscr{C})$
			is bounded above. Then, for any pullback square
		\[
				\begin{tikzcd}
						P' \arrow[d] \arrow[r]\arrow[dr, phantom, "\lrcorner", very near start] & P \arrow[d] \\
						Q' \arrow[r]           & Q         
				\end{tikzcd}
		\]
		in $CAlg^{pd,ele}_{B//(A/I)}$ for which the maps $\pi_0(P)\rightarrow \pi_0(Q)$ and $\pi_0(Q')\rightarrow\pi_0(Q)$ are surjective, applying $Def_M^{pd}$ gives a pullback square in $\hat{S}$. The same assertion holds for the sheaf-algebra deformation functor $Def_X$.
	\end{proposition}
	\begin{proof}
			Suspension induces an equivalence between $Def^{pd}_M$ and $Def^{pd}_{M[d]}$. Thus, after a shift, we may assume that $M$ is connective.
			\\By Proposition~\ref{conservative}, every $N\in Def^{pd}_M(P)$ is connective, so $Def^{pd}_M(P)$ is a subcategory of $Mod^{cn}_{P}(\mathscr{C})$.
			\\By \cite[Theorem 16.2.0.2]{SAG}, the map $Mod^{cn}_{P'}(\mathscr{C})\rightarrow Mod^{cn}_{P}(\mathscr{C})\times_{Mod^{cn}_{Q}(\mathscr{C})}Mod^{cn}_{Q'}(\mathscr{C})$ is an equivalence.
			\\Consequently, the commutative diagram
		\[
		\begin{tikzcd}
			Def^{pd}_M(P') \arrow[d] \arrow[r] & Def^{pd}_M(P) \arrow[d] \\
			Def^{pd}_M(Q') \arrow[r]           & Def^{pd}_M(Q)          
		\end{tikzcd}
		\]
		is a pullback square.
			\par For $Def_X$, we glue a compatible pair of sheaf-algebra deformations by taking their pullback as sheaves of animated rings. The same patching argument for the underlying connective sheaves of modules shows that completed base change recovers the original pair and that a deformation over $P'$ is recovered from its two base changes.
			\par The forgetful functor detects equivalences and preserves limits, which proves the sheaf-algebra assertion as in the proof of \cite[Proposition 3.25]{CY}. For a split square-zero extension of $A/I$, the square-zero fiber sequence gives the required connectivity and conservativity directly, so the same argument applies to these test objects even when they are outside the stated elementary subcategory.
	\end{proof}
		The deformation functors and Proposition~\ref{rigidity} give the following commutative diagram for every $C\in CAlg^{pd,ele}_{B//(A/I)}$, as in \cite[Chapter 6.(15)]{CY}. The same construction applies to the split square-zero objects in Proposition~\ref{rigidity}:
	\begin{equation}\label{per0}
\begin{tikzcd}
	Def_X(C) \arrow[r, "\bullet\hat{\otimes}^L_{und(C)} C"] & Def^{pd}_X(C) \arrow[d] \arrow[rrr, "{\varphi^*R\Gamma_{\Prism}(X)\hat{\otimes}^L_A C\rightarrow R\Gamma(X,\bullet)}"] &  &  &  Def_{R\Gamma(\theta_X)}(C) \arrow[d] \\
	& \{*\}=\{X\} \arrow[rrr, "dR^{pd}_{X/C}"]                                                                           & &  & Def_{R\Gamma(dR^{pd}_{X/(A/I)})}(C) 
\end{tikzcd}
	\end{equation} 
		Here $\theta_X$ denotes the natural map $dR^{pd}_{X/(A/I)}\rightarrow O^{pd}_X$, and the map $\varphi^*R\Gamma_{\Prism}(X)\hat{\otimes}^L_A C\rightarrow R\Gamma(X,\bullet)$ is induced by $\varphi^*R\Gamma_{\Prism}(X)\hat{\otimes}^L_A B\simeq R\Gamma_{cris}(X)$ and Proposition~\ref{rigidity}.
	
		We now define the period map, which sends a deformation of the sheaf $O_X$ to a deformation of a morphism of complexes.
		\begin{definition}(\cite[Definition 6.13]{CY})
			Define $Per_X$ by the following pullback square:
		\begin{equation}
			\begin{tikzcd}
				Per_X \arrow[d] \arrow[r] \arrow[dr, phantom, "\lrcorner", very near start] &  Def_{R\Gamma(\theta_X)} \arrow[d]   \\
				\{*\} \arrow[r]           & Def_{R\Gamma(dR^{pd}_{X/(A/I)})}
			\end{tikzcd}
		\end{equation}
			The period map $per_X: Def_X\rightarrow Per_X$ is the map induced by (\ref{per0}).
	\end{definition}
		We will use a slightly different map to relax some of the freeness assumptions in the proof below.
	\begin{definition}
			Define $Per^1_X$ by the following pullback square:
		\begin{equation}
			\begin{tikzcd}
			Per^1_X \arrow[d] \arrow[r] \arrow[dr, phantom, "\lrcorner", very near start] &  Def_{R\Gamma(\theta^1_X)} \arrow[d]   \\
			\{*\} \arrow[r]           & Def_{R\Gamma(dR^{pd}_{X/(A/I)})}
			\end{tikzcd}
		\end{equation}
			Here $R\Gamma(\theta^1_X)$ is the composite map $R\Gamma(dR^{pd}_{X/(A/I)})\rightarrow R\Gamma(O_X)\rightarrow H^d(X,O_X)[-d]$.
	\end{definition}
		\begin{proposition}
			Let $X$ be a formal scheme satisfying Assumption~\ref{ass}, and let $C\in CAlg^{pd,ele}_{B//(A/I)}$ be discrete. There is a canonical map $f:Def_X(C)\rightarrow Per^1_X(C)$, given by sending $O_{\tilde{X}}\in Def^{pd}_X(C)$ to $H^d(X,O_{\tilde{X}}^{pd})$.
	\end{proposition}
	\begin{proof}
			By Assumption~\ref{ass}(a), the prismatic cohomology admits a decomposition
		\[
		\varphi^*R\Gamma_{\Prism}(X)\cong A[-d]\oplus M
		\]
		for some $M\in\widehat D(A)$, since a retract of a derived $p$-complete complex is derived $p$-complete.
			\\Grothendieck's vanishing theorem gives
			\[	\tau_{\leqslant -d}R\Gamma(X,O_{\tilde{X}}^{pd})[d]\cong H^d(X,O_{\tilde{X}}^{pd})\in D^{\heartsuit}(C).
		\]
			By hypothesis, $R\Gamma(X,O_{\tilde{X}}^{pd})$ and the discrete ring $C$ are derived $p$-complete. Hence $H^d(X,O_{\tilde{X}}^{pd})$ and the kernels and cokernels of maps between these complete modules are derived $p$-complete.
			\\Since completed derived tensor products preserve cohomological upper bounds and $R\Gamma(X,O_{\tilde{X}}^{pd})\in D^{\leqslant d}(C)$, for any discrete $C$-module $T$ we have
		\[
		H^d(X,O_{\tilde{X}}^{pd}\hat{\otimes}_C^L T)\simeq H^0(H^d(X,O_{\tilde{X}}^{pd})\hat{\otimes}_C^L T)=H^d(X,O_{\tilde{X}}^{pd})\hat{\otimes}_C T.
		\]
			In particular, the base-change identification associated with the deformation, $R\Gamma(X,O_{\tilde{X}}^{pd})\hat{\otimes}_C^L A/I\simeq R\Gamma(X,O_X)$, gives $H^d(X,O_{\tilde{X}}^{pd})\hat{\otimes}_C A/I\simeq H^d(X,O_X)\simeq A/I$.
			\\A diagram chase therefore reduces the claim to showing that the map
		\[
		\varphi^*R\Gamma_{\Prism}(X)\hat{\otimes}_A^L C\rightarrow \tau_{\leqslant -d}R\Gamma(X,O_{\tilde{X}}^{pd})
		\] lies in $Def_{R\Gamma(\theta_X^1)}$.
			\\Equivalently, we must show that
		\begin{equation}\label{basechange}
			H^d(X,O_{\tilde{X}}^{pd})[-d]\hat{\otimes}^L_C A/I\cong H^d(X,O_X)[-d].
		\end{equation}
			By Assumption~\ref{ass}(a), we may write $\varphi^*R\Gamma_{\Prism}(X)=A[-d]\oplus M$.
			\\By definition, $B\rightarrow C$ is surjective, so the kernel of $C\rightarrow A/I$ is $JC$.
			\\Consider the following commutative diagram:
\begin{equation}\label{O_X}
	\begin{tikzcd}
				(\tau_{\geqslant 1-d}M)\hat{\otimes}_A^L (JC) \arrow[d] \arrow[r] & (\tau_{\geqslant 1-d}M)\hat{\otimes}_A^L C \arrow[r] \arrow[d] & (\tau_{\geqslant 1-d}M)\hat{\otimes}_A^L (A/I) \arrow[d, "h"] \\
		R\Gamma(X,O_{\tilde{X}}^{pd})\hat{\otimes}_C^L(JC) \arrow[r]        & R\Gamma(X,O_{\tilde{X}}^{pd}) \arrow[r]                          & {R\Gamma(X,O_X)}     
	\end{tikzcd}
\end{equation}
			The horizontal rows are exact.
			\\By Assumption~\ref{ass}(b), $H^{d-1}(h)$ is surjective. Consider the induced cohomology diagram:
		\[
\begin{tikzcd}
	H^{d-1}((\tau_{\geqslant 1-d}M)\hat{\otimes}_A^L (A/I)) \arrow[d, "H^{d-1}(h)"] \arrow[r] & H^d((\tau_{\geqslant 1-d}M)\hat{\otimes}_A^L (JC))=0 \arrow[d] \\
	{H^{d-1}(X,O_X)} \arrow[r]                                                                & {H^d(R\Gamma(O_{\tilde{X}}^{pd})\hat{\otimes}_C^L(JC))}                                                
\end{tikzcd}
		\]
			The horizontal maps are induced by (\ref{O_X}). It follows that the boundary map
		\[
		H^{d-1}(X,O_X)\rightarrow H^d(R\Gamma(O_{\tilde{X}}^{pd})\hat{\otimes}_C^L(JC))= H^d(X,O_{\tilde{X}}^{pd})\hat{\otimes}_C(JC)
		\]
		is zero.
			\\We therefore obtain an injection
		\[
			0\rightarrow H^d(X,O_{\tilde{X}}^{pd})\hat{\otimes}_C(JC)\rightarrow H^d(X,O_{\tilde{X}}^{pd}).
		\]
			\\Now consider the map $C\rightarrow H^d(X,O_{\tilde{X}}^{pd})$ induced by
		\[
		C[-d]\rightarrow (C[-d]\oplus M\hat{\otimes}_A^L C)\cong \varphi^*R\Gamma_{\Prism}(X)\hat{\otimes}_A^L C\rightarrow H^d(X,O_{\tilde{X}}^{pd})[-d]
		\]
		and denote its cokernel in $\widehat D^{\heartsuit}(C)$ by $N$.
			\\After base change to $A/I$, we obtain the exact sequence
		\[
			A/I=C\hat{\otimes}_C A/I\rightarrow A/I=H^d(X,O_X)=H^d(X,O_{\tilde{X}}^{pd})\hat{\otimes}_C A/I\rightarrow N\hat{\otimes}_C A/I\rightarrow 0.
		\]
			The first map is the identity under the chosen identification with $A/I$. Hence $H^0(N\hat{\otimes}_C^L A/I)=0$. Applying Proposition~\ref{conservative} to $N[-1]$ gives $N=0$, so $C\rightarrow H^d(X,O_{\tilde{X}}^{pd})$ is surjective. Its kernel $J_1\in\widehat D^{\heartsuit}(C)$ lies in $JC$.
		
		The injection $H^d(X,O_{\tilde{X}}^{pd})\hat{\otimes}_C(JC)\hookrightarrow H^d(X,O_{\tilde{X}}^{pd})$ above, together with the triangle obtained by tensoring $JC\rightarrow C\rightarrow A/I$ with $H^d(X,O_{\tilde{X}}^{pd})$, implies $H^{-1}(H^d(X,O_{\tilde{X}}^{pd})\hat{\otimes}_C^L A/I)=0$. Tensoring $J_1\rightarrow C\rightarrow H^d(X,O_{\tilde{X}}^{pd})$ with $A/I$ now gives an exact sequence
		\[
		0\longrightarrow H^0(J_1\hat{\otimes}_C^L A/I)\longrightarrow A/I\xrightarrow{\sim}H^0(H^d(X,O_{\tilde{X}}^{pd})\hat{\otimes}_C^L A/I)\simeq A/I.
		\]
			Thus $H^0(J_1\hat{\otimes}_C^L A/I)=0$. Applying Proposition~\ref{conservative} to $J_1[-1]$ gives $J_1=0$.
			\\Therefore, $H^d(X,O_{\tilde{X}}^{pd})\simeq C$, which establishes the full derived base-change identification (\ref{basechange}).
	\end{proof}
		To prove our main theorem, we need an explicit description of the period map, particularly on tangent spaces over $A/I$.
	\begin{proposition}
			For a smooth proper $p$-adic formal scheme $X$ over $A/I$, any discrete derived $p$-complete $A/I$-module $V$, and $n\geqslant0$, the period map $per_X(sqz^{pd}_{A/I}(V[n]))$ fits into the following commutative square:
		\begin{equation}\label{percalculate}
					\begin{tikzcd}
				{Def_X(sqz^{pd}_{A/I}(V[n]))} \arrow[d] \arrow[r, "\simeq"] & {Map_{O_X}(L_{X/(A/I)},V[n+1]\hat{\otimes}^L O_X)} \arrow[d]             \\
				{Per_X(sqz^{pd}_{A/I}(V[n]))} \arrow[r, "\simeq"]           & {Map_{A/I}(R\Gamma(X,Fil^1dR_{X/(A/I)}),V[n]\hat{\otimes}^L R\Gamma(X,O_X))}
			\end{tikzcd}
		\end{equation}
			Here $L_{X/(A/I)}$ denotes the derived $p$-completed cotangent complex; for smooth $X$, it is the finite locally free sheaf of continuous differentials in degree zero. Set $Fil^1dR_{X/(A/I)}=Fib(dR_{X/(A/I)}\rightarrow O_X)$. The right vertical map is induced by
		\[
			\pi^1:Fil^1dR_{X/(A/I)}\rightarrow L_{X/(A/I)}[-1].
		\]
	\end{proposition}
	\begin{proof}
			Write $C:=sqz^{pd}_{A/I}(V[n])$ and $N:=R\Gamma(dR^{pd}_{X/(A/I)})$.
			\\First, Proposition~\ref{formalmoduli} gives
		\[
			Def_{N}(C)=\Omega Def_{N}(C[1]).
		\]
			Unraveling the definition gives the following pullback square:
		\[
\begin{tikzcd}
	\Omega Def_{N}(C[1]) \arrow[r] \arrow[d]\arrow[dr, phantom, "\lrcorner", very near start] & {Map_{Mod^{\simeq}_{C[1]}}(C[1]\hat{\otimes}^L_{A/I}N,C[1]\hat{\otimes}^L_{A/I}N)} \arrow[d, "A/I\hat{\otimes}^L_{C[1]}\bullet"] \\
	\{*\} \arrow[r, "id"]                 & {Map_{Mod_{A/I}^{\simeq}}(N,N)}                
\end{tikzcd}
		\]
			Here $Mod^{\simeq}_{A/I}$ denotes the maximal $\infty$-groupoid of the category $Mod^{\wedge}_{A/I}$: it contains all objects of $Mod^{\wedge}_{A/I}$ and only the equivalences between $A/I$-modules.
			\\For any $N'\in Mod^{\wedge}_{C[1]}$, we also have a fiber sequence
		\[
			V[n+1]\hat{\otimes}^L_{A/I}(A/I\hat{\otimes}^L_{C[1]}N')\rightarrow N'\rightarrow A/I\hat{\otimes}^L_{C[1]}N'.
		\]
			Thus the base-change functor $N'\mapsto A/I\hat{\otimes}^L_{C[1]}N'$ is conservative, giving the following pullback square:
		\[
		\begin{tikzcd}
			{Map_{Mod^{\simeq}_{C[1]}}({C[1]}\hat{\otimes}^L_{A/I}N,{C[1]}\hat{\otimes}^L_{A/I}N)} \arrow[d, "A/I\hat{\otimes}^L_{C[1]}\bullet"] \arrow[r] \arrow[dr, phantom, "\lrcorner", very near start]& {Map_{Mod_{C[1]}}({C[1]}\hat{\otimes}^L_{A/I}N,{C[1]}\hat{\otimes}^L_{A/I}N)} \arrow[d, "A/I\hat{\otimes}^L_{C[1]}\bullet"] \\
			{Map_{Mod_{A/I}^{\simeq}}(N,N)} \arrow[r]                                & {Map_{Mod_{A/I}}(N,N)}                                                                                 
		\end{tikzcd}
		\]
			Combining this with the equivalence $Map_{Mod_C}(C\hat{\otimes}^L_{A/I}N,C\hat{\otimes}^L_{A/I}N)\cong Map_{Mod_{A/I}}(N,C\hat{\otimes}^L_{A/I}N)$
			gives
		\[
			Def_{N}(C)\cong Map'_{Mod_{A/I}}(N,C[1]\hat{\otimes}^L_{A/I} N).
		\]
			For a specified morphism $y\to x$, let $Map'(x,y)$ denote the homotopy fiber over $id_x$ of the induced map $Map(x,y)\to Map(x,x)$, as in the following pullback square:
		\[
		\begin{tikzcd}
			{Map'(x,y)} \arrow[d] \arrow[r]\arrow[dr, phantom, "\lrcorner", very near start] & {Map(x,y)} \arrow[d] \\
			\{*\} \arrow[r, "id"]           & {Map(x,x)}          
		\end{tikzcd}
		\]
		Similarly, we have
		\[
			Def_{R\Gamma(\theta_X)}(C)\cong Map'_{Mod_{A/I}}(R\Gamma(\theta_X),C[1]\hat{\otimes}^L_{A/I} R\Gamma(\theta_X)),
		\]
		\[
			Def_X(C)=Map'_{{Shv(\mathscr{T},CAlg^{pd})}}(O_X, O_X\hat{\otimes}^L_{A/I}C[1]).
		\]
			By definition, the commutative square
		\[
	\begin{tikzcd}
		Def^{pd}_X(C) \arrow[d] \arrow[r] &  \{*\}\arrow[d] \\
		Def_{R\Gamma(\theta_X)}(C) \arrow[r]                   & Def_{R\Gamma(dR^{pd}_{X/(A/I)})}(C) 
	\end{tikzcd}
		\]
			 in (\ref{per0}) is equivalent to the outer square of the following diagram:
		 \[
		\begin{tikzcd}
			{Map'(O_X, O_X\hat{\otimes}^L_{A/I}C[1])} \arrow[r] \arrow[d, "R\Gamma\circ\theta"]                          & \{*\} \arrow[d]                                                        \\
			{Map'(R\Gamma(\theta_X), R\Gamma(\theta_{X\hat{\otimes}^L_{A/I}C[1]}) )} \arrow[d] \arrow[r]        & \{*\} \arrow[d]                                                        \\
			{Map'(R\Gamma(\theta_X),C[1]\hat{\otimes}^L_{A/I} R\Gamma(\theta_X))} \arrow[r] & {Map'(R\Gamma(dR_{X/(A/I)}),C[1]\hat{\otimes}^LR\Gamma(dR_{X/(A/I)}))}
		\end{tikzcd}
		 \]
			 Here the lower $\{*\}$ corresponds to $Map'(R\Gamma(dR_{X/(A/I)}),R\Gamma(dR_{X/(A/I)}))$, and $R\Gamma(\theta_{X\hat{\otimes}^L_{A/I}C[1]})$ denotes the complex $\left[  R\Gamma(dR_{X/(A/I)})\rightarrow R\Gamma(O_X)\hat{\otimes}^L_{A/I}C[1]\ \right]$.
			 \\Comparing fibers shows that the lower square is a pullback square.
			 \\Hence (\ref{percalculate}) is equivalent to the upper square of this diagram.
			 \\Finally, we have the following two equivalences:
		 \[
			 (1)Map'(O_X, O_X\hat{\otimes}^L_{A/I}C[1])\simeq Map_{O_X}(L_{X/(A/I)},V[n+1]\hat{\otimes}^L O_X),
		 \]
		 \[
		 \begin{split}
		 	(2)&Map'(R\Gamma(\theta_X),\left[  R\Gamma(dR_{X/(A/I)})\rightarrow R\Gamma(O_X)\hat{\otimes}^L_{A/I}C[1]\ \right]  )
		 	\\&\simeq Map(R\Gamma(\theta_X),\left[  0\rightarrow R\Gamma(O_X)\hat{\otimes}^L_{A/I}V[n+1]\ \right]  )
			 \\&\simeq Map_{A/I}(R\Gamma(X,Fil^1dR_{X/(A/I)}),V[n]\hat{\otimes}^L R\Gamma(X,O_X)).
		 \end{split}
		 \]
			 The first equivalence follows from the properties of the cotangent complex for square-zero extensions. The second follows from the decomposition
		 \[
		 \begin{gathered}
	\left[  R\Gamma(dR_{X/(A/I)})\rightarrow R\Gamma(O_X)\hat{\otimes}^L_{A/I}C[1]\ \right] =\\
		 \left[  R\Gamma(dR_{X/(A/I)})\rightarrow R\Gamma(O_X)\ \right] \oplus \left[  0\rightarrow R\Gamma(O_X)\hat{\otimes}^L_{A/I}V[n+1]\ \right]
		 \end{gathered}
		 \]
		 and the fiber sequence
		 \[
			 R\Gamma(X,Fil^1dR_{X/(A/I)})\rightarrow R\Gamma(dR_{X/(A/I)})\rightarrow R\Gamma(O_X).
		 \]
		 \\Thus, the corresponding map
		 \[
		 Map_{O_X}(L_{X/(A/I)},V[n+1]\hat{\otimes}^L O_X)\rightarrow Map_{A/I}(R\Gamma(X,Fil^1dR_{X/(A/I)}),V[n]\hat{\otimes}^L R\Gamma(X,O_X))
		 \]
			 is induced by the commutative square
		 \[
		 \begin{tikzcd}
		 	dR_{X/(A/I)} \arrow[r] \arrow[d] & O_X \arrow[d, "d"] \\
		 	0 \arrow[r]                      & L_{X/(A/I)}       
		 \end{tikzcd}
		 \]
			 or, equivalently, by the map $\pi^1:Fil^1dR_{X/(A/I)}\rightarrow L_{X/(A/I)}[-1]$. This proves the proposition.
	\end{proof}
	\begin{proposition}\label{per1}
			Let $X$ be a formal scheme satisfying Assumption~\ref{ass}, and let $C$ be an object of $CAlg^{pd,ele}_{B//(A/I)}$. Write $\varphi^*R\Gamma_{\Prism}(X)=A[-d]\oplus M$. Then $\pi_0Per^1_X(C)\simeq Hom_A(M,I_C[-d])$.
			\\Here $I_C=Fib(C\rightarrow A/I)$ is viewed as an $A$-module via $A\rightarrow C$.
	\end{proposition}
	\begin{proof}
			By Proposition~\ref{conservative}, for any $\varphi^*R\Gamma_{\Prism}(X)\hat{\otimes}^L_A C\rightarrow N$ deforming $R\Gamma_{dR}(X)\rightarrow H^d(X,O_X)$, the composite map
		\[
		C[-d]\rightarrow C[-d]\oplus M\hat{\otimes}^L_A C\rightarrow N
		\]
		is an equivalence.
			\\Let $Mod^{\simeq,*}_{C,C\hat{\otimes}^L\varphi^*R\Gamma_{\Prism}(X)/}$ be the subcategory of $Mod^{\wedge}_C$ whose objects are maps $C[-d]\oplus M\hat{\otimes}^L_A C\rightarrow Q$ for which the induced map $C[-d]\rightarrow Q$ is an equivalence, and whose morphisms are equivalences between such objects.
			\\Unraveling the definition, we see that $Per^1_X(C)$ fits into the following pullback square:
		\[
		\begin{tikzcd}
			Per^1_X(C) \arrow[d] \arrow[r]\arrow[dr, phantom, "\lrcorner", very near start] & {Mod^{\simeq,*}_{C,C\hat{\otimes}^L\varphi^*R\Gamma_{\Prism}(X)/}} \arrow[d] \\
			\{*\} \arrow[r]                & {Mod^{\simeq,*}_{A/I,R\Gamma_{dR}(X)}}    
		\end{tikzcd}
		\]
			Moreover, we have
		\[
			Mod^{\simeq,*}_{C,C\hat{\otimes}^L\varphi^*R\Gamma_{\Prism}(X)/}\simeq Map_C(M\hat{\otimes}^L_A C, C[-d])\simeq Map_A(M,C[-d]).
		\]
		Similarly,
		\[
			Mod^{\simeq,*}_{A/I,R\Gamma_{dR}(X)}\simeq Map_A(M,A/I[-d]).
		\]
			We therefore obtain the fiber sequence
		\[
			Per^1_X(C)\rightarrow Map_A(M,C[-d])\rightarrow Map_A(M,A/I[-d]).
		\]
		Finally, since $Map_A(M,-)$ preserves fibers, the preceding fiber sequence gives
		\[
		Per^1_X(C)\simeq Map_A(M,I_C[-d]).\qedhere
		\]

	\end{proof}
	Fix a decomposition $\varphi^*R\Gamma_{\Prism}(X)=A[-d]\oplus M$ whose reduction is compatible with the top Hodge quotient. Set
	\[
	M_{A/I}:=M\hat{\otimes}^L_A R\simeq Fib\bigl(R\Gamma(dR_{X/R})\longrightarrow H^d(X,O_X)[-d]\bigr).
	\]
		For a discrete $p$-complete $A/I$-module $V$, write
		\[
		\begin{aligned}
		M_1:=Hom_{A/I}(&R\Gamma(X,Fil^1dR_{X/(A/I)}),\\
		&H^d(X,O_X)[-d]\hat{\otimes}^L_{A/I} V[1]).
		\end{aligned}
		\]
	
		For $ sqz^{pd}_{A/I}(V[1])$, we have the following maps:
	\[
	\begin{split}
		f_1: \pi_0Per_X(sqz^{pd}_{A/I}(V[1]))\simeq Hom_{A/I}(R\Gamma(X,Fil^1dR_{X/(A/I)}),V[1]\hat{\otimes}^L_{A/I} R\Gamma(X,O_X))\rightarrow M_1,
	\end{split}
	\]
	\[
	\begin{split}
		f_2:\pi_0Per^1_X(sqz^{pd}_{A/I}(V[1]))&\simeq Hom_A(M,V[1-d])\\
			&\simeq Hom_{A/I}(M_{A/I},V[1-d])\rightarrow M_1.
	\end{split}
	\]
		Here $f_1$ is induced by $R\Gamma(X,O_X)\rightarrow H^d(X,O_X)[-d]$. Under the displayed adjunction, $f_2$ is given by precomposition with the natural map $R\Gamma(X,Fil^1dR_{X/(A/I)})\rightarrow M_{A/I}$, using the identification $H^d(X,O_X)\simeq A/I$ fixed above.
		\\These two maps fit into the following commutative diagram.
	\begin{proposition}\label{per1map}
			Let $C\in CAlg^{pd,ele}_{B//(A/I)}$ be discrete, and let $C\rightarrow sqz^{pd}_{A/I}(V[1])$ be a map, where $V$ is a discrete $p$-complete $A/I$-module. Then the following diagram commutes:
		\[
		\begin{tikzcd}
			\pi_0Def_X(C) \arrow[r] \arrow[d] & {\pi_0Def^{pd}_X(sqz^{pd}_{A/I}(V[1]))} \arrow[r]   & {\pi_0Per_X(sqz^{pd}_{A/I}(V[1]))} \arrow[d,"f_1"]                                \\
			\pi_0Per_X^1(C) \arrow[r]         & {\pi_0Per_X^1(sqz^{pd}_{A/I}(V[1]))} \arrow[r,"f_2"] & M_1
		\end{tikzcd}
		\]
			Here $M_1:=Hom_{A/I}(R\Gamma(X,Fil^1dR_{X/(A/I)}),H^d(X,O_X)[-d]\hat{\otimes}^L_{A/I} V[1])$, as above.
	\end{proposition}
	\begin{proof}
			Write $D$ for $sqz^{pd}_{A/I}(V[1])$.
			For any $O_{\tilde{X}}\in\pi_0Def_X(C)$, we compare its images in
		\[
Hom_{A/I}(R\Gamma(X,Fil^1dR_{X/(A/I)}),H^d(X,O_X)[-d]\hat{\otimes}^L_{A/I} V[1])
			\] along the two paths.
			\\Along the upper path, it maps to
		\[
		R\Gamma(dR_{X/D})\rightarrow R\Gamma(O^{pd}_{\tilde{X}}\hat{\otimes}^L_C D)
		\]
			in $\pi_0Per_X(D)$.
			The commutative diagram of fiber sequences
		\[
		\begin{tikzcd}
			{R\Gamma(Fil^1dR_{X/(A/I)}\hat{\otimes}^L_{A/I} D)} \arrow[r] \arrow[d] & {R\Gamma(dR_{X/D})} \arrow[d] \arrow[r]                    & {R\Gamma(O^{pd}_X)\hat{\otimes}^L_{A/I} D} \arrow[d] \\
			{R\Gamma(O^{pd}_X)\hat{\otimes}^L_{A/I} V[1]} \arrow[r]         & {R\Gamma(O^{pd}_{\tilde{X}}\hat{\otimes}^L_C D)} \arrow[r] & R\Gamma(O^{pd}_X)                                                      
		\end{tikzcd}
		\]
		induces, by precomposition with the section $A/I\rightarrow D$, the map
		\[
		R\Gamma(Fil^1dR_{X/(A/I)})\rightarrow R\Gamma(O^{pd}_X)\hat{\otimes}^L_{A/I} V[1]
		\]
			which represents the same element under the equivalence
		\[
			\pi_0Per_X(D)\simeq Hom_{A/I}(R\Gamma(X,Fil^1dR_{X/(A/I)}),V[1]\hat{\otimes}^L_{A/I} R\Gamma(X,O_X)).
		\]
			Composing with $f_1$ yields the map
		\begin{equation}\label{prsqz}
				R\Gamma(Fil^1dR_{X/(A/I)})\rightarrow R\Gamma(O^{pd}_X)\hat{\otimes}^L_{A/I} V[1]\rightarrow H^d(X,O_X)[-d]\hat{\otimes}^L_{A/I} V[1].
		\end{equation}
			\\Along the lower path, $O_{\tilde{X}}$ maps to
		\[
		R\Gamma(dR_{X/D})\rightarrow H^d(X,O^{pd}_{\tilde{X}})[-d]\hat{\otimes}^L_C D
		\]
			in $\pi_0Per_X^1(D)$.
		Using (\ref{basechange}), the same fiber construction and precomposition with $A/I\rightarrow D$ give the map
		\[
			M_{A/I}\rightarrow H^d(X,O_X)[-d]\hat{\otimes}^L_{A/I} V[1].
		\]
			\\Composing with $f_2$ yields the map
		\[
		R\Gamma(Fil^1dR_{X/(A/I)})\rightarrow M_{A/I}\rightarrow H^d(X,O_X)[-d]\hat{\otimes}^L_{A/I} V[1]
		\]
		which agrees with (\ref{prsqz}) by functoriality of fibers. Indeed, the natural quotient $R\Gamma(X,O^{pd}_{\tilde{X}})\rightarrow H^d(X,O^{pd}_{\tilde{X}})[-d]$, after base change to $D$ and reduction to $A/I$, induces the map between augmentation fibers given by $R\Gamma(X,O_X)\rightarrow H^d(X,O_X)[-d]$ tensored with $V[1]$. On the source fibers it induces $R\Gamma(X,Fil^1dR_{X/(A/I)})\rightarrow M_{A/I}$. Thus the two composites define the same element of $M_1$.
			\\Thus the diagram in the statement of the proposition commutes.
		\end{proof}
	\section{The existence of deformations}
	In this section, we prove Theorem A using the period map constructed above. We assume that Assumption \ref{ass} holds and use the chosen compatible decomposition $\varphi^*R\Gamma_{\Prism}(X)=A[-d]\oplus M$.
	\begin{proposition}\label{inj}
		Assume that $(A,I)$ is a transversal or crystalline prism. If $X$ satisfies Assumption \ref{ass} and its Hodge--de Rham spectral sequence degenerates at $E_1$, then the map
		\begin{equation}\label{period}
			\pi_0Def_X(sqz^{pd}_{A/I}(I^n/I^{n+1}[1]))\rightarrow Hom(R\Gamma(X,Fil^1dR_{X/{(A/I)}}),H^d(X,O_X)[-d]\hat{\otimes} F_n/F_{n+1}[1])
		\end{equation} induced by the period map is injective. Here $(B,J,\gamma)$ is the $p$-complete PD-envelope of $(A,I)$, and $F_n=\overline{J^n}$ is the closure of $J^n$ in $B$.
	\end{proposition}
	\begin{proof}
		First, we prove that the map $I^n/I^{n+1}\rightarrow F_n/F_{n+1}$ is split.
		\\When $(A,I)$ is crystalline, the claim follows immediately from $I^n=J^n=F_n$.
		\\When $(A,I)$ is transversal, let $D_0$ be the uncompleted PD-envelope of $(A,I)$, so that $B=\widehat{D_0}_p$.
		\\There is a canonical map
		\[
		D_0\longrightarrow A/I^m[\frac{1}{p}],\qquad \gamma_k(a)\longmapsto\frac{a^k}{k!}\quad(a\in I).
		\]
	This map can be extend uniquely to a ring homomorphism
		\[
		\phi_m:B\longrightarrow A/I^m[\frac{1}{p}].
		\]
		These maps are compatible with reduction in $m$, and $\phi_1$ is the augmentation $B\rightarrow R$ followed by $R\hookrightarrow R[1/p]$.
		
		For $m=2$, the image lies in $A/I^2$ itself, and $\phi_2(J)\subset I/I^2$. For each $x\in J$, choose $a_x\in I$ lifting $\phi_2(x)$. Compatibility gives
		\[
		\phi_{n+1}(x)-a_x\in I^2A/I^{n+1}[\frac{1}{p}].
		\]
		Thus, for $x_1,\ldots,x_n\in J$, we have
		\[
		\phi_{n+1}(x_1\cdots x_n)=a_{x_1}\cdots a_{x_n}\in I^n/I^{n+1}\subset A/I^{n+1}[\frac{1}{p}],
		\]
		and $\phi_{n+1}(J)\subset IA/I^{n+1}[\frac{1}{p}]$, so $\phi_{n+1}(J^{n+1})=0$.
		\\The submodule $I^n/I^{n+1}$ is closed in $A/I^{n+1}[\frac{1}{p}]$; continuity therefore gives $\phi_{n+1}(F_n)\subset I^n/I^{n+1}$ and $\phi_{n+1}(F_{n+1})=0$. Since $JF_n\subset F_{n+1}$, restriction induces a canonical $R$-linear map
		\[
		s_n:F_n/F_{n+1}\longrightarrow I^n/I^{n+1}.
		\]
		This splits the map $I^n/I^{n+1}\rightarrow F_n/F_{n+1}$. Thus, $I^n/I^{n+1}$ is a direct summand of $F_n/F_{n+1}$.
		
		Unwinding the definition, we obtain the following description of the map (\ref{period}):
		\begin{equation}
			\begin{split}
			&Hom_{O_X}(\Omega^1[-1],O_X\{n\}[1])\\&\simeq  Hom(R\Gamma(X,\Omega^1_X),H^d(X,O_X)[2-d])\{n\}\\&\stackrel{f_1}{\rightarrow} Hom(R\Gamma(X,Fil^1dR_{X/{(A/I)}}),H^d(X,O_X)[1-d])\{n\}\\&\stackrel{f_2}{\rightarrow} H^0(RHom(R\Gamma(X,Fil^1dR_{X/{(A/I)}}),H^d(X,O_X)[-d]\otimes F_n/F_{n+1}[1]))
			\end{split}
		\end{equation}
		Here the first equivalence follows from Serre duality. The map $f_1$ is induced by $Fil^1dR_{X/(A/I)}\rightarrow \Omega^1_{X/(A/I)}[-1]$. The map $f_2$ is induced by the natural injection $I^n/I^{n+1}\hookrightarrow F_n/F_{n+1}$.
		\\Moreover, Serre duality and Poincaré duality give
		\[
		Hom(R\Gamma(X,Fil^1dR_{X/{(A/I)}}),H^d(X,O_X)[1-d])\simeq H^{d+1}(X,\Omega^{<d}_{X/(A/I)})
		\]
		\[
		Hom(R\Gamma(X,\Omega^1_{X/(A/I)}),H^d(X,O_X)[2-d])\simeq H^{2}(X,\Omega^{d-1}_{X/(A/I)})
		\]
		By the degeneration of the Hodge--de Rham spectral sequence, the map
		\[
		H^{2}(X,\Omega^{d-1}_{X/(A/I)})\rightarrow H^{d+1}(X,\Omega^{<d}_{X/(A/I)})
		\]
		is injective. Tensoring with the invertible module $I^n/I^{n+1}$ preserves injectivity, so $f_1$ is injective.
		\\Finally, $f_2$ is injective because $I^n/I^{n+1}$ is a direct summand of $F_n/F_{n+1}$.
	\end{proof}

	\begin{proposition}\label{surj}
		Let $(A,I)$ be a transversal or crystalline prism, and let $m$ be the Tor-dimension of $A$ ($m$ can be $\infty$). Let $X$ be a proper smooth $p$-adic formal scheme over $A/I$ of relative dimension $d$, with geometrically connected fibres and trivial relative canonical bundle. Assume the following:
		\\(a) $H^i_{\Prism}(X)$ is finite projective for $i\in[d,d+m]$;
		\\(b) The Hodge--de Rham spectral sequence of $X$ degenerates at $E_1$.
		\\Then the natural map $\pi_0Def_X(A/I^{n+1})\rightarrow \pi_0Def_X(A/I^n)$ is surjective for any $n>1$. It is also surjective for $n=1$ if $(A,I)$ is transversal, or if $(A,I)$ is crystalline with $p>2$.
	\end{proposition}
	\begin{proof}
		Since $X$ is proper and smooth with geometrically connected fibres, $H^0(X,O_X)=A/I$, relative Serre duality then gives $H^d(X,O_X)\cong A/I$.
		\\By \cite[Theorem 6.4]{duality},
		\[
		(R\Gamma_{\Prism}(X))^\vee\cong R\Gamma_{\Prism}(X)[2d]\{d\}.
		\]
		Here the twist $\{d\}$ is an invertible $A$-module by \cite[Notation 2.5.4]{absolute}.
		\\Since $H^i_{\Prism}(X)$ is finite projective for $i\in[d,d+m]$, the Tor spectral sequence gives
		\[
		\begin{aligned}
		H^d(\varphi^*R\Gamma_{\Prism}(X)\otimes_A^L N)&\cong H^d_{\Prism}(X)\otimes_{A,\varphi}N,\\
		Hom_A(\varphi^*R\Gamma_{\Prism}(X),N[-d])&\cong H^d((R\Gamma_{\Prism}(X))^\vee[-2d])\otimes_{A,\varphi}N.
		\end{aligned}
		\]
		for any $A$-module $N$.
		\\The second equivalence shows that $Hom_A(\varphi^*R\Gamma_{\Prism}(X),A/I[-d])$ preserves surjections of discrete coefficient modules.
		\\Hence, the top Hodge quotient $q:\varphi^*R\Gamma_{\Prism}(X)\to A/I[-d]$ lifts to a map $\lambda:\varphi^*R\Gamma_{\Prism}(X)\to A[-d]$, which gives a splitting $\varphi^*R\Gamma_{\Prism}(X)\cong A[-d]\oplus M$ and proves Assumption \ref{ass}(a).
		\\To verify Assumption \ref{ass}(b), note that the Cartier divisor $I$ gives $A/I$ Tor-amplitude $[-1,0]$ over $A$. Since $H^d(\varphi^*R\Gamma_{\Prism}(X))$ is projective, tensoring the cohomological truncation triangle shows that
		\[
		(\tau_{\geqslant1-d}\varphi^*R\Gamma_{\Prism}(X))\otimes_A^L A/I\cong\tau_{\geqslant1-d}(\varphi^*R\Gamma_{\Prism}(X)\otimes_A^L A/I)=\tau_{\geqslant1-d}R\Gamma_{dR}(X/(A/I)).
		\]
		Hodge--de Rham degeneration makes the map from the right-hand side to $\tau_{\geqslant1-d}R\Gamma(X,O_X)$ surjective in every degree. The displayed complexes are derived $p$-complete, so the same identity holds with completed tensor products. Thus Assumption \ref{ass}(b) holds as well.
		\par First, suppose $n=1$ and $(A,I)$ is transversal. Equip $A/I^2$ with the PD structure of Lemma \ref{firstthickening}, and let $V=I/I^2$. Proposition \ref{per1} identifies
		\[
		\pi_0Per_X^1(A/I^2)=Hom_A(M,V[-d]),\qquad Per_X^1(A/I)\simeq\{*\}.
		\]
		Thus the period datum of $X$ lifts to $A/I^2$ via the zero element on the left, or equivalently via the projection $\varphi^*R\Gamma_{\Prism}(X)\hat\otimes_A^L A/I^2\to A/I^2[-d]$.
		\par The pullback square of Lemma \ref{firstthickening} defines the obstruction
		\[
		o_X\in\pi_0Def_X(sqz_{A/I}^{pd}(V[1]))\cong H^2(X,T_{X/(A/I)}\otimes_{A/I} V)
		\]
		to lifting $X$ to $A/I^2$. Apply Proposition \ref{per1map} to the discrete base $A/I$ and its obstruction map $A/I\to sqz_{A/I}^{pd}(V[1])$. It identifies the image of $o_X$ under the period map with the image of the period obstruction. The latter vanishes because the period datum lifts. The map
		\[
		H^2(X,T_{X/(A/I)}\otimes_{A/I} V)\longrightarrow
		Hom_{A/I}(R\Gamma(X,Fil^1dR_{X/(A/I)}),V[1-d])
		\]
		is injective by the $f_1$ argument in the proof of Proposition \ref{inj}. Hence $o_X=0$, proving the assertion for $n=1$ in the transversal case.
		\par We now treat $n\geqslant2$, together with $n=1$ in the crystalline case with $p>2$. Let $(B,J,\gamma)$ be the $p$-complete PD-envelope of $(A,I)$. Since $I^nB\subset F_n$, there are compatible maps of ordinary rings $A/I^n\to B/F_n$. Proposition \ref{pullback} supplies the PD elementary extension in these cases. Using the period map, we obtain the following commutative diagram with exact rows:
		\begin{equation}\label{exact}
			\begin{tikzcd}
				\pi_0 Def_X(A/I^{n+1}) \arrow[r] \arrow[d] & \pi_0 Def_X(A/I^n) \arrow[r] \arrow[d] & {\pi_0Def_X(sqz^{pd}_{A/I}(I^n/I^{n+1}[1]))} \arrow[d] \\
				\pi_0 Def_X(B/F_{n+1}) \arrow[r] \arrow[d] & \pi_0 Def_X(B/F_n) \arrow[r] \arrow[d] & {\pi_0Def_X(sqz^{pd}_{B/J}(F_n/F_{n+1}[1]))} \\
				\pi_0 Per^1_X(B/F_{n+1}) \arrow[r]         & \pi_0 Per^1_X(B/F_n) \arrow[r]         & {\pi_0Per^1_X(sqz^{pd}_{B/J}(F_n/F_{n+1}[1]))}        
			\end{tikzcd}
		\end{equation}
		By Propositions \ref{per1map} and \ref{inj}, it therefore suffices to show that $\pi_0 Per^1_X(B/F_{n+1}) \rightarrow \pi_0 Per^1_X(B/F_n)$ is surjective.
		\\By Proposition \ref{per1}, this map corresponds to the following:
		\[
		Hom_A(M,(J/F_{n+1})[-d])\rightarrow Hom_A(M,(J/F_n)[-d])
		\]
		The functor $Hom_A(M,-[-d])$ is a natural direct summand of $Hom_A(K,-[-d])$. The latter preserves surjections by the Tor calculation above, so the former does as well. Applying this to $J/F_{n+1}\twoheadrightarrow J/F_n$ proves the required surjectivity. Equivalently, taking the corresponding direct summand of the natural tensor description gives
		\[
		Hom_A(M,N[-d])\cong H^{-d}(M^\vee)\otimes_A N=H^d(M^\vee[-2d])\otimes_A N.
		\]
		Here $H^{-d}(M^\vee)$ is finite projective, and for the complete coefficient modules $N=J/F_r$, the tensor product on the right is already derived $p$-complete.
	\end{proof}
	\Needspace{12\baselineskip}
	\begin{theorem}\label{mainlift}
		Let $A$ be an integral ring, and let $(A,I)$ be a transversal prism or a crystalline prism with $p>2$. Let $m$ be the Tor-dimension of $A$. If $X$ is a smooth proper $p$-adic formal scheme over $A/I$ of relative dimension $d$ with trivial canonical bundle and satisfies the following conditions:
		\begin{center}
				(a) $H^i_{\Prism}(X)$ is free for all $i\in[d,d+m]$;
				\\(b) The Hodge--de Rham spectral sequence of $X$ degenerates at $E_1$,
		\end{center}
		then there exists a compatible system of smooth proper $p$-adic formal deformations $X_n$ over $A/I^n$, for all $n\geqslant1$, with $X_1\cong X$. This system defines a smooth proper formal deformation $\mathscr{X}$ of $X$ over $\operatorname{Spf}_{(p,I)}A$.
	\end{theorem}
	\begin{proof}
		We first reduce to geometrically connected fibres. Put $S=H^0(X,O_X)$. The Stein factorization makes $S$ finite \'etale over $A/I$, and $X\to\operatorname{Spf}(S)$ has geometrically connected fibres.
		
		The finite \'etale algebra $S$ lifts uniquely to a finite \emph{\'etale} $A$-algebra $A'$ with $A'/IA'=S$, and the $\delta$-structure extends uniquely, giving a transversal or crystalline prism $(A',IA')$, which does not increase the Tor-dimension of $A$.
		
		Therefore, without loss of generality, we may assume that $X$ has geometrically connected fibres, with finite projective cohomology in the required range.
		
		We inductively construct a compatible system of sheaves $\mathscr{F}_n\in Shv(X,CAlg^{an})$ over $A/I^n$, starting with $\mathscr{F}_1=O_X$.
		
		Suppose $\mathscr{F}_n$ has been constructed, with $n\geqslant1$. Applying Proposition \ref{surj} directly to $(A,I)$, including its first-step assertion when $n=1$, we choose a lift $\mathscr{F}_{n+1}\in\pi_0Def_X(A/I^{n+1})$ together with an identification
		\[
		\mathscr{F}_{n+1}\hat{\otimes}^L_{A/I^{n+1}}A/I^n\cong\mathscr{F}_n.
		\]
		\\By induction, we have a compatible system of sheaves $\mathscr{F}_n\in Shv(X, CAlg^{an})$ over $A/I^n$ such that $\mathscr{F}_{n+1}\hat{\otimes}^L_{A/I^{n+1}}A/I^{n}=\mathscr{F}_n$.
		\\In particular, $\pi_0\mathscr{F}_{n+1}\hat{\otimes}_{A/I^{n+1}}A/I^{n}=\pi_0\mathscr{F}_n$.
		\\By Remark \ref{classicalthickenings}, each $\mathscr{F}_n$ is discrete and represents a smooth proper $p$-adic formal scheme over $A/I^n$. More precisely, for an affine open $U=\operatorname{Spf}R_0\subset X$, the ring $R_{n,U}=\Gamma(U,\pi_0\mathscr{F}_n)$ is $p$-complete, and the restriction of $\pi_0\mathscr{F}_n$ to $U$ is $O_{\operatorname{Spf}_p R_{n,U}}$. The restrictions of the given sheaves supply the gluing and the compatibility between successive levels.
		\par Thus $X_n=(|X|,\pi_0\mathscr{F}_n)$ is the required compatible system of smooth proper $p$-adic formal deformations over $A/I^n$. Put
		\[
		\mathscr{X}=\left(|X|,\varprojlim_n\pi_0\mathscr{F}_n\right),
		\]
		where the structure sheaf is endowed with the $(p,I)$-adic topology. To justify the formal-scheme assertion, put $\mathfrak a=(p)+I$ and let $Y_m$ be the base change of $X_m$ to $A/\mathfrak a^m$. These form a compatible system of smooth proper schemes. On an affine open $U\subset X$, the ideals $(p^r,I^n)$ are cofinal with powers of $\mathfrak a$, so
		\[
		R_U:=\varprojlim_n R_{n,U}\cong\varprojlim_m\Gamma(U,O_{Y_m}).
		\]
		Since $\mathfrak a$ is finitely generated, \cite[\href{https://stacks.math.columbia.edu/tag/09B8}{09B8}]{stack} gives $R_U/\mathfrak a^mR_U\cong\Gamma(U,O_{Y_m})$. Restrictions to distinguished opens are the corresponding $\mathfrak a$-completed localizations, identifying $(U,\mathcal{O}_{\mathscr{X}}|_U)$ with $\operatorname{Spf}_{(p,I)}R_U$. These charts glue by the restrictions of the given sheaves. Smoothness and properness follow from the smooth proper reductions $Y_m$. Consequently, $\mathscr{X}$ is a smooth proper formal scheme over $\operatorname{Spf}_{(p,I)}A$, whose base change to $\operatorname{Spf}_p(A/I^n)$ is $X_n$ for every $n\geqslant1$.
	\end{proof}
	\Needspace{6\baselineskip}
	\begin{remark}[Crystalline prisms at $p=2$]\label{char2}
		For a crystalline prism $(A,(2))$, the conclusion of Theorem \ref{mainlift} still holds if an initial deformation over $A/4$ exists: Proposition \ref{surj} then applies at every subsequent step. The first-step PD argument does not apply, since the canonical divided powers satisfy $\gamma_2(2)=2\not\equiv0\pmod4$.
		\par If $A=W(k)$ with $k$ a perfect field of any positive characteristic, Hodge--de Rham degeneration supplies the initial lift by \cite[Theorem 7.18]{CY}, after the same reduction to geometrically connected fibres. Thus no additional hypothesis is needed over $W(k)$ at any prime. In particular, if $k$ is a perfect field of characteristic $p>0$ and $X$ is a smooth proper $k$-scheme of pure dimension $d$ with trivial canonical bundle and satisfies the following conditions:
		\begin{center}
			(a) $H^i_{crys}(X/W(k))$ is free for $i=d,d+1$;
			\\(b) The Hodge--de Rham spectral sequence of $X$ degenerates at $E_1$,
		\end{center}
		then $X$ admits a smooth proper $p$-adic formal lift over $\operatorname{Spf}_pW(k)$.
		\\These two conditions are weaker than those in \cite[Theorem A]{CY}, which requires all crystalline cohomology groups to be free. They are precisely conditions (NCT) and (DEG) in \cite{ST}.
	\end{remark}
	\begin{example}[Complete intersections in Totaro's fivefold]\label{totaroexample}
		Let $k$ be an algebraically closed field of characteristic $2$,. On $G=\operatorname{Gr}(2,k^4)$, let $\mathcal S$ be the tautological rank-two bundle, and put
		\[
		\pi:V=\mathbf P_G(F_G^*\mathcal S)\longrightarrow G,
		\qquad L=\mathcal O_V(1)\otimes\pi^*\mathcal O_G(1),
		\]
		where projectivization parametrizes lines and $F_G$ is the absolute Frobenius. This is the fivefold in \cite[\S1]{Totaro}, with $(n,p)=(4,2)$. The line bundle $L$ is very ample and $\omega_V\simeq L^{-2}$. Thus the intersection $Y$ of two general members of the complete linear system $|L|$ is a smooth projective threefold with $\omega_Y\simeq\mathcal O_Y$.
		\par Projective-bundle and Koszul calculations, together with the normal sequence, give
		\[
		h^q(Y,\mathcal O_Y)=1,0,0,1,
		\qquad h^q(Y,T_Y)=0,32,2,0,\ \text{for}\ q=0,1,2,3.
		\]
		Here generality is understood in the full space $H^0(V,L)$, which has dimension $25$.
		\\On the other hand, for any $l\not=2$
		\[
		b_n(Y)=\dim_{\mathbb {Q}_l}H^n_{\mathrm{\acute et}}(Y,\mathbb {Q}_l)=1,0,2,66,2,0,1, \text{for $n$ ranging from 0 to 6}
		\]
		 The identifications $\Omega_Y^2\simeq T_Y$ and $\Omega_Y^3\simeq\mathcal O_Y$, together with Serre duality, show that the total Hodge dimensions equal the Betti numbers in every degree. Crystalline base change and the Hodge--de Rham spectral sequence therefore give
		\[
		b_n(Y)\leqslant\dim_kH^n_{dR}(Y/k)
		\leqslant\sum_{a+b=n}h^b(Y,\Omega_Y^a)=b_n(Y).
		\]
		Equality forces $E_1$-degeneration and freeness of all crystalline cohomology groups; in particular,
		\[
		H^3_{crys}(Y/W(k))\simeq W(k)^{66},\qquad H^4_{crys}(Y/W(k))\simeq W(k)^2.
		\]
		\par The degeneration also supplies an initial lift to $W_2(k)$ by \cite[Theorem 7.18]{CY}. Remark \ref{char2} now applies to the crystalline prism $(W(k),(2))$, with $d=3$ and $m=1$: freeness in degrees $3$ and $4$ ensures that the refined period data lift, and injectivity of the period map on obstructions then yields compatible smooth proper lifts of $Y$ over all $W_n(k)$. These define a formal lift over $\operatorname{Spf}W$. Finally, $H^2(Y,\mathcal O_Y)=0$ allows an ample line bundle to lift through this tower, so the formal lift algebraizes to a smooth projective $W$-scheme.
	\end{example}
	\section{Applications to \texorpdfstring{$F$}{F}-split Calabi--Yau varieties}
	We now apply our construction to deformations of $F$-split varieties over an algebraically closed field $k$ of characteristic $p>0$. Throughout this section, we fix the prism $(A,I)$ to be $(W(k),p)$. Petrov's theorem \cite[Theorem 1.1]{Pet25} implies that the Hodge--de Rham spectral sequence of every smooth proper $F$-split variety considered here degenerates at $E_1$, without any dimension restriction. We use this degeneration below without imposing it as a separate hypothesis. We begin with the following proposition for Calabi--Yau varieties.
	\begin{proposition}\cite[Proposition 2.5.1]{ST}\label{st}
		Let $X$ be a smooth, proper variety of dimension $d$ over $k$. Assume that $X$ has trivial canonical bundle and that $H^{d-1}(X,O_X)=0$. Then the following are equivalent:
		\\(1) $X$ admits a Frobenius splitting.
		\\(2) $F:H^d(X,WO_X)\rightarrow H^d(X,WO_X) $ is bijective.
		\\Moreover, under these equivalent conditions, the natural map
		\[
		H^d(X,WO_X)^{F=1}\otimes_{\mathbb{Z}_p} W(k)\rightarrow H^d(X,WO_X)
		\]
		is an isomorphism, and $H^d(X,WO_X)$ is free of rank one over $W(k)$.
	\end{proposition}
	\begin{proof}
		The Frobenius map gives a morphism
		\[
		\iota: O_X\rightarrow F_*O_X,
		\]
		which induces a map
		\[
		\iota': Hom_{O_X}(F_*O_X,O_X)\rightarrow Hom_{O_X}(O_X,O_X)\cong k.
		\]
		By Serre duality, this map is dual to the map
		\[
		F: H^d(X,O_X)\rightarrow H^d(X,F_*O_X).
		\]
		Thus $X$ admits a Frobenius splitting $\Leftrightarrow$ $\iota'$ is bijective $\Leftrightarrow$ $F$ is bijective.
		\\If these conditions hold, the filtration and the vanishing $H^{d-1}(X,O_X)=0$ show by induction that $H^{d-1}(X,W_nO_X)=0$ for every $n\geqslant 0$.
		\\Together with the exact sequence
		\[
		0\rightarrow W_nO_X\xrightarrow{V} W_{n+1}O_X\rightarrow O_X\rightarrow 0,
		\]
		this gives the exact sequence
		\[
		0\rightarrow H^{d}(X, W_nO_X)\xrightarrow{V} H^d(X,W_{n+1}O_X)\rightarrow k\rightarrow 0.
		\]
		\\Since $F$ is bijective and $p=FV$, the module $H^d(X,WO_X)$ is $p$-torsion-free and satisfies
		\[
		H^d(X,WO_X)/pH^d(X,WO_X)=k.
		\]
		\\Therefore, $H^d(X,WO_X)$ is free of rank one.
		\\Finally, since $F$ is bijective and $k$ is algebraically closed, induction modulo $p$ gives an element $x$ in $H^d(X,WO_X)$ satisfying $F(x)=x$, which yields an $F$-fixed basis.
	\end{proof}
	\begin{proposition}(cf.\ \cite[Proposition 2.4.1]{ST})\label{st2}
		Let $X$ be an $F$-split smooth, proper variety of dimension $d$ over $k$. Assume that $X$ has trivial canonical bundle and that $H^{d-1}(X,O_X)=0$. Then there is an $F$-stable decomposition
		\[
		H^d_{cris}(X/W(k))=H^d(X,WO_X)\oplus H_{\geqslant 1},
		\]
		satisfying $F(H_{\geqslant 1})\subset pH_{\geqslant 1}$.
		\\Here $F: H^d_{cris}(X/W(k))\rightarrow H^d_{cris}(X/W(k))$ is the natural Frobenius map on crystalline cohomology.
	\end{proposition}
	\begin{proof}
		 The degree-zero projection of the de Rham--Witt complex induces a Frobenius-equivariant map $\pi:H^d_{cris}(X/W(k))\rightarrow H^d(X,WO_X)$. By Proposition \ref{st}, $H^d(X,WO_X)$ is free of rank one, its Frobenius is bijective, and $H^d(X,WO_X)/pH^d(X,WO_X)\cong H^d(X,O_X)=k$.
		\par The compatible Frobenius Fitting decompositions give
		\[
		H^d_{cris}(X/W(k))=U\oplus T,
		\]
		where $F$ is bijective on $U$ and topologically nilpotent on $T$. Frobenius equivariance implies $\pi(T)=0$: for every $r$, a sufficiently high iterate sends $T$ into $p^rT$, whereas the Frobenius map preserves $p^rH^d(X,WO_X)$ bijectively and $H^d(X,WO_X)$ is separated. The rational slope comparison \[H^d(X,WO_X)[1/p]\cong H^d_{cris}(X/W(k))[1/p]_{<1}
		\]
		in \cite[Introduction, (1.1)]{BBE} shows that $U$ has rank one, since $H^d(X,WO_X)[1/p]$ has slope zero and $U[1/p]$ is the slope-zero part of $H^d_{cris}(X/W(k))[1/p]$.
		\par Crystalline base change gives a Frobenius-equivariant injection
		\[
		\alpha:H^d_{cris}(X/W(k))/pH^d_{cris}(X/W(k))\hookrightarrow H^d_{dR}(X/k).
		\]
		On the de Rham complex, absolute Frobenius factors as
		\[
		\Omega^\bullet_{X/k}\xrightarrow{\mathrm{pr}_0}O_X[0]
		\xrightarrow{a\mapsto a^p}F_{X*}\Omega^\bullet_{X/k},
		\]
		Write $\rho:H^d_{dR}(X/k)\rightarrow H^d(X,O_X)$ for the degree-zero projection. Compatibility of the projections gives $\rho(\alpha(\bar x))=\pi(x)\bmod p$ for $x\in H^d_{cris}(X/W(k))$.
		\\If $u\in U$ satisfies $\pi(u)\in pH^d(X,WO_X)$, this factorization and injectivity of $\alpha$ imply $F(u)\in pH^d_{cris}(X/W(k))\cap U=pU$. Since $F$ is bijective on $U$, we have $u\in pU$. Thus $\pi$ induces an injection
		\[
		U/pU\hookrightarrow H^d(X,WO_X)/pH^d(X,WO_X)=k.
		\]
		As $U$ has rank one, the structure theorem for finite $W(k)$-modules shows that $U$ is free of rank one and this injection is an isomorphism. Nakayama's lemma now gives a Frobenius-equivariant isomorphism $\pi|_U:U\xrightarrow{\sim}H^d(X,WO_X)$. Hence $H_{\geqslant 1}:=\ker\pi=T$ is an $F$-stable complement.
		\par Finally, for $x\in\ker\pi$, the same de Rham factorization gives $F(x)\in pH^d_{cris}(X/W(k))$. Write $F(x)=py$ with $y\in H^d_{cris}(X/W(k))$. Then
		\[
		p\pi(y)=\pi(F(x))=0.
		\]
		Since $H^d(X,WO_X)$ is $p$-torsion-free, $\pi(y)=0$, so $y\in\ker\pi$. This proves $F(H_{\geqslant 1})\subset pH_{\geqslant 1}$.
	\end{proof}
	\begin{remark}
		Proposition \ref{st2} gives the conclusion of \cite[Proposition 2.4.1]{ST} for $m=1$ under the present hypotheses. The proof above does not require the crystalline-freeness condition (NCT) or Hodge--de Rham degeneration: the divisibility argument uses only the injection $H/pH\hookrightarrow H^d_{dR}(X/k)$ and the torsion-freeness of $H^d(X,WO_X)$.
	\end{remark}
	This proposition yields the following corollary.
	\begin{corollary}
		Let $X$ be an $F$-split smooth, proper variety of dimension $d$ over $k$. Assume that $X$ has trivial canonical bundle and that $H^{d-1}(X,O_X)=0$. Then $X$ satisfies Assumption \ref{ass}.
	\end{corollary}
	\begin{proof}
		By \cite[Theorem 1.1]{Pet25}, the Hodge--de Rham spectral sequence of $X$ degenerates at $E_1$. The de Rham--Witt complex of $X$ induces a map
	\[
	R\Gamma_{cris}(X/W(k))\rightarrow R\Gamma(X,WO_X)\rightarrow R\Gamma(X,O_X).
	\]
	Since $H^i(X,WO_X)=0$ for $i>d$, this induces a map
	\[
	R\Gamma_{cris}(X/W(k))\rightarrow H^d(X,WO_X)[-d]\cong W(k)[-d]\rightarrow H^d(X,O_X)[-d].
	\]
	On the other hand, since $W(k)$ has dimension one, we have an exact sequence
	\begin{equation}\label{spectral}
		0 \rightarrow H^i_{cris}(X/W(k))\otimes_{W(k)}k \rightarrow H^i_{dR}(X) \rightarrow {Tor^1_{W(k)}(H^{i+1}_{cris}(X/W(k)),k)} \rightarrow 0.
	\end{equation}
	Together with Hodge--de Rham degeneration, this shows that, to verify Assumption \ref{ass}(b), it suffices to check degree $d-1$.
	\\Let $M$ be the cokernel of $H^{d-1}_{cris}(X)\hat{\otimes}^L_{W(k)} k\rightarrow H^{d-1}(X,O_X)$. It remains to show that $M=0$.
	\\By (\ref{spectral}) and Hodge--de Rham degeneration, we have a surjective map
	\begin{equation}\label{Tormap}
		Tor^1_{W(k)}(H^d_{cris}(X/W(k)),k)\rightarrow M.
	\end{equation}
	By Proposition \ref{st2}, we can decompose $Tor^1_{W(k)}(H^d_{cris}(X/W(k)),k)$ as
	\[
	Tor^1_{W(k)}(H^d(X,WO_X),k)\oplus Tor^1_{W(k)}(H_{\geqslant 1},k).
	\]
	The first summand vanishes because $H^d(X,WO_X)$ is free. For the second summand, compatibility with Frobenius gives a commutative diagram
	\[
	\begin{tikzcd}
		{Tor^1(H_{\geqslant1},k)} \arrow[d, "F^n"] \arrow[r] & M \arrow[d, "F^n"] \\
		{Tor^1(H_{\geqslant1},k)} \arrow[r]             &  M           
	\end{tikzcd}
	\]
	for every $n>0$, where $F^n$ denotes the $n$-fold iterate of Frobenius.
	\\The crystalline cohomology of $X$ is finitely generated. Together with $F(H_{\geqslant1})\subset pH_{\geqslant1}$, this gives $p^n(H_{\geqslant1})_{tor}=0$ for sufficiently large $n$. Thus $F^n(H_{\geqslant1})\subseteq p^nH_{\geqslant1}$ is torsion-free.
	\\For such $n$, the induced map $F^n$ on $Tor^1(H_{\geqslant1},k)$ is therefore zero, since this is a $p$-torsion module.
	\\On the other hand, the existence of a Frobenius splitting $\sigma$ shows that the Frobenius map on $H^{d-1}(X,O_X)$ is bijective and hence induces a bijective map $F: M\rightarrow M$.
	\\It follows that the upper horizontal map in the diagram is zero. Thus the map (\ref{Tormap}) is zero, so $M=0$.
	\\Together with Hodge--de Rham degeneration, this shows that $X$ satisfies Assumption \ref{ass}.
	\end{proof}
	In \cite[Section 5]{ST}, the authors study the following class of deformations of $F$-split varieties.
	\begin{definition}\label{canonical}
		Let $X$ be an $F$-split smooth, proper Calabi--Yau variety over $k$, and let $\tilde{X}_n$ be a deformation of $X$ over $W_n(k)$. Define the Hodge filtration $F^1H^d_{dR}(\tilde{X}_n/W_n(k))$ as the image of the map
		\[
		H^d(\tilde{X}_n,\Omega^{\bullet\geqslant1})\rightarrow H^d_{dR}(\tilde{X}_n/W_n(k)).
		\]
		We call such a deformation \textit{Hodge–Frobenius compatible} if $F^1H^d_{dR}(\tilde{X}_n/W_n(k))$ is Frobenius-stable. Here the Frobenius map on $H^d_{dR}(\tilde{X}_n/W_n(k))$ is defined via the canonical map
		\[
		H^d_{dR}(\tilde{X}_n/W_n(k))\cong H^d_{cris}(X/W_n(k)).
		\]
	\end{definition}
	\begin{remark}
		For $X$ as in Definition \ref{canonical}, assume that $p>2$ and that $H^d_{cris}(X/W(k))$ and $H^{d+1}_{cris}(X/W(k))$ are free. Hodge--de Rham degeneration follows from \cite[Theorem 1.1]{Pet25}. Then, for each $n\geqslant1$, there is at most one Hodge–Frobenius compatible deformation over $W_n(k)$ up to isomorphism, by \cite[Theorem 5.7.1, Remark 5.7.3]{ST}. For ordinary abelian varieties and ordinary K3 surfaces, Hodge–Frobenius compatible deformations exist for every $n$ and are the reductions of their classical canonical lifts; see \cite[Corollary 5.7.2, Remark 5.7.3]{ST}.
	\end{remark}
	Under the hypotheses below, our construction shows that any specified Hodge--Frobenius compatible lift over $W_n(k)$ extends to $W_{n+1}(k)$.
	\begin{theorem}\label{lift}
		 Let $X$ be an $F$-split smooth, proper variety of dimension $d$ over the field $k$ of characteristic $p>0$. Assume that $X$ has trivial canonical bundle and that $H^{d-1}(X,O_X)=0$. If a Hodge–Frobenius compatible lift $\tilde{X}_n$ exists over $W_n(k)$ for $n\geqslant1$, then $\tilde{X}_n$ admits a lift to $W_{n+1}(k)$.
	\end{theorem}
	\begin{proof}
		By \cite[Theorem 1.1]{Pet25}, the Hodge--de Rham spectral sequence of $X$ degenerates at $E_1$, so the degeneration hypothesis required in the obstruction comparison is satisfied. For $n=1$, the Frobenius-splitting construction $W_2(O_X)/I_\sigma$ gives a lift of $X$ over $W_2(k)$ in every characteristic, by \cite[Theorems 3.16 and 3.18]{zda}. Thus we may assume $n\geqslant2$, when Proposition \ref{pullback} supplies the required PD elementary-extension square for every prime $p$.
		\par\Needspace{7\baselineskip}
		As in (\ref{exact}), we have the following commutative diagram:
		\begin{equation}
			\begin{tikzcd}
				\pi_0Def_X(W_{n+1}(k)) \arrow[r] \arrow[d] & \pi_0Def_X(W_{n}(k)) \arrow[r] \arrow[d] & {\pi_0Def_X(sqz^{pd}_k(k[1]))}  \\
				\pi_0Per^1_X(W_{n+1}(k)) \arrow[r]         & \pi_0Per^1_X(W_{n}(k))\arrow[r]          & {\pi_0Per^1_X(sqz^{pd}_k(k[1]))}        
			\end{tikzcd}
		\end{equation}
		with the horizontal row exact.
		\\Similar with Proposition \ref{surj}, we can show that if the image of $\pi_0Def_X(W_{n}(k)$ in $\pi_0Per^1_X(sqz^{pd}_k(k[1]))$ is zero, then the image of $\pi_0Def_X(W_{n}(k)$ in $\pi_0Def_X(sqz^{pd}_k(k[1]))$ is also zero.
		\\Thus, to show that $\tilde{X}_n\in \pi_0Def_X(W_n(k))$ lies in the image of $\pi_0Def_X(W_{n+1}(k))$, it suffices to show that the map
		\begin{equation}\label{6.2}
			R\Gamma_{cris}(X/W(k))\otimes^L_{W(k)}W_n(k)\rightarrow H^d(X,O_{\tilde{X}_n})[-d]
		\end{equation}
		between complexes over $W_n(k)$ lifts to $W_{n+1}(k)$.
		\\On the Zariski site, the exact sequence of abelian sheaves
		\[
		0\longrightarrow WO_X\xrightarrow{V^n}WO_X\longrightarrow W_nO_X\longrightarrow0
		\]
		and the vanishing $H^{d+1}(X,WO_X)=0$ give
		\[
		H^d(X,WO_X)/p^nH^d(X,WO_X)=H^d(X,WO_X)/V^nH^d(X,WO_X)\xrightarrow{\sim}H^d(X,W_nO_X).
		\]
		These identifications commute with Frobenius and restriction in $n$.
		\\Derived crystalline base change \cite[II, Theorem 2.7 and Scholie 2.8]{Ill79} gives the Frobenius-equivariant exact sequence
		\[
		0\longrightarrow H^d(X/W(k))/p^nH^d(X/W(k))\xrightarrow{i_n}H^d_{cris}(X/W_n(k))
		\longrightarrow H^{d+1}(X/W(k))[p^n]\longrightarrow0
		\]
		Since for sufficiently large $e$, $\varphi^e(H^{d+1}(X/W(k))[p^n])=0$, and the bijectivity of Frobenius on $H^d(X,WO_X)/p^nH^d(X,WO_X)$, we can view $H^d(X,W_nO_X)$ as a direct summand of $H^d_{cris}(X/W_n(k))$.
		\\Now, we show that $H^{d}(\tilde{X}_n,O_{\tilde{X}_n})$ is free of rank one. For the reductions $\tilde{X}_j$ of $\tilde{X}_n$, the successive structure-sheaf sequences and the vanishing $H^{d-1}(X,O_X)=0$ give, by induction,
		\[
		H^{d-1}(\tilde{X}_j,O_{\tilde{X}_j})=0,\qquad
		0\longrightarrow k\longrightarrow H^{d}(\tilde{X}_j,O_{\tilde{X}_j})\longrightarrow H^{d}(\tilde{X}_{j-1},O_{\tilde{X}_{j-1}})\longrightarrow0
		\quad(2\leqslant j\leqslant n).
		\]
		Thus $H^{d}(\tilde{X}_n,O_{\tilde{X}_n})$ has length $n$. If $u\in H^{d}(\tilde{X}_n,O_{\tilde{X}_n})$ lifts a generator of $H^{d}(X,O_{X})=k$, then $p^{n-1}u\ne0$ by the injection $k\hookrightarrow H^{d}(\tilde{X}_n,O_{\tilde{X}_n})$ above. Hence $u$ freely generates $H^{d}(\tilde{X}_n,O_{\tilde{X}_n})$ over $W_n(k)$.
		\\Now denote $\pi_n: H^d_{cris}(X/W_n(k))\rightarrow H^d(X,W_nO_X)$ and $q_n:H^d_{cris}(X/W_n(k))\to H^{d}(\tilde{X}_n,O_{\tilde{X}_n})$ as the edge map. 
		\\The restriction $q_n|_{H^d(X,W_nO_X)}$ reduces modulo $p$ to the isomorphism \[H^d(X,WO_X)/pH^d(X,WO_X)\to H^d(X,O_X)\] established in Proposition \ref{st2}.
		\\Both $H^d(X,W_nO_X)$ and $H^{d}(\tilde{X}_n,O_{\tilde{X}_n})$ are free of rank one. Hence $q_n|_{H^d(X,W_nO_X)}$ is an isomorphism. In particular, $q_n$ is surjective and
		\[
		H^d(X,W_nO_X)\cap F^1H^d_{dR}(\tilde{X}_n/W_n(k))
		=H^d(X,W_nO_X)\cap\ker q_n=0.
		\]
		Now write an element of this Hodge filtration as $x=x_1+x_2$, with $x_1\in H^d(X,W_nO_X)$ and $x_2\in \ker \pi_n$. For $m\geqslant n+e$, Hodge–Frobenius compatibility gives
		\[
		\varphi^m(x)=\varphi^m(x_1)\in H^d(X,W_nO_X)\cap\ker q_n=0.
		\]
		Bijectivity of $\varphi$ on $H^d(X,W_nO_X)$ implies $x_1=0$. Hence $\ker q_n\subset N_n$. The corresponding quotients of $H^d_{cris}(X/W_n(k))$ are $H^{d}(\tilde{X}_n,O_{\tilde{X}_n})$ and $H^d(X,W_nO_X)$, both of length $n$, so this inclusion is an equality. We obtain
		\[
		F^1H^d_{dR}(\tilde{X}_n/W_n(k))=\ker\pi_n
		\]
		and a canonical isomorphism
		\[
		\lambda_n:H^d(X,W_nO_X)\xrightarrow{\sim}H^{d}(\tilde{X}_n,O_{\tilde{X}_n}),\qquad q_n=\lambda_n\circ\pi_n.
		\]
		Here $\lambda_n=q_n\circ j_n$ under the identification $H^d(X,W_nO_X)\cong U/p^nU$.
		The isomorphism $\lambda_n$ therefore identifies the canonical map
		\[
		R\Gamma_{cris}(X/W_n(k))\rightarrow R\Gamma(X,W_n(O_X))\rightarrow H^d(X,W_n(O_X))[-d]\rightarrow H^d(X,O_{\tilde{X}_n})[-d]
		\]
		with (\ref{6.2}). Indeed, their maps on $H^d$ agree by $q_n=\lambda_n\circ\pi_n$, and $H^{d}(\tilde{X}_n,O_{\tilde{X}_n})$ is injective over $W_n(k)$, so equality can be checked on $H^d$.
		\\Therefore, the map
		\[
		R\Gamma_{cris}(X/W_{n+1}(k))\rightarrow R\Gamma(X,W_{n+1}(O_X))\rightarrow H^d(X,W_{n+1}(O_X))[-d]
		\]
		has the free rank-one target $H^d(X,W_{n+1}O_X)[-d]$. Witt restriction identifies its reduction with $H^{d}(\tilde{X}_n,O_{\tilde{X}_n})[-d]$ via $\lambda_n$. Compatibility of crystalline base change and Witt restriction shows that it reduces to (\ref{6.2}); thus it is a deformation of that map over $W_{n+1}(k)$.
		\\This shows that the image of $\tilde{X}_n$ in $\pi_0Per^1_X(sqz^{pd}_k(k[1]))$ is zero. Hence $\tilde{X}_n$ lifts to $W_{n+1}(k)$.
	\end{proof}
	In particular, we obtain the following liftability result over $W_3(k)$.
	\begin{corollary}(Theorem B)
		Let $X$ be an $F$-split smooth proper variety of dimension $d$ over an algebraically closed field $k$ of characteristic $p>2$, with trivial canonical bundle and $H^{d-1}(X,O_X)=0$. Then $X$ admits a lift to $W_3(k)$.
	\end{corollary}
	\begin{proof}
		Choose an $F$-splitting $\sigma:F_*O_X\to O_X$ and set $O_{\tilde{X}}=W_2(O_X)/I_\sigma$, where $I_\sigma=\{(0,f)\mid\sigma(f)=0\}$. This construction commutes with localization and, on every affine open, gives a flat $W_2(k)$-algebra lifting $O_X$, by \cite[Theorem 3.16, Lemma 3.17, Theorem 3.18]{zda}. The affine lifts therefore glue to a scheme $\tilde{X}=(|X|,O_{\tilde{X}})$ lifting $X$. Smoothness and properness follow from those of the special fibre across the nilpotent thickening \cite[\href{https://stacks.math.columbia.edu/tag/06AG}{06AG}]{stack}.
		\par To verify Hodge–Frobenius compatibility for proper $X$, set $K=\ker\sigma$. The splitting gives $F_*O_X\cong O_X\oplus K$. Frobenius is injective on each finite-dimensional $H^q(X,O_X)$ and hence bijective, since $k$ is perfect. Thus $H^q(X,K)=0$ for every $q$. As additive sheaves, $I_\sigma\cong K$, so the quotient induces isomorphisms
		\[
		\beta_q:H^q(X,W_2(O_X))\xrightarrow{\sim}H^q(\tilde{X},O_{\tilde{X}}),
		\]
		as in \cite[Lemma 4.2.3]{ST}.
		\par Equip the square-zero ideal $VO_X=\ker(W_2(O_X)\to O_X)$ with trivial divided powers. Since $p>2$, these extend the canonical divided powers on $(p)\subset W_2(k)$, whose higher divided powers vanish modulo $p^2$. The quotient $W_2(O_X)\to O_{\tilde{X}}$ is then a morphism of PD thickenings of $X$. Crystalline restriction, followed by projection to degree zero, gives
		\[
		H^d_{cris}(X/W_2(k))\xrightarrow{\alpha}H^d(X,W_2(O_X))
		\xrightarrow{\beta_d}H^d(\tilde{X},O_{\tilde{X}});
		\]
		see \cite[\href{https://stacks.math.columbia.edu/tag/07LP}{07LP}]{stack}. By functoriality and crystalline--de Rham comparison, the composite is the Hodge edge map. The Witt-vector Frobenius is a morphism of the first PD thickening over the Witt Frobenius of $W_2(k)$, so $\alpha$ is Frobenius-equivariant. Since $\beta_d$ is an isomorphism, it follows that
		\[
		F^1H^d_{dR}(\tilde{X}/W_2(k))=\ker(\beta_d\circ\alpha)=\ker\alpha
		\]
		is Frobenius-stable. Thus $\tilde{X}$ is Hodge–Frobenius compatible in the sense of Definition \ref{canonical}, and Theorem \ref{lift} with $n=2$ gives the required lift to $W_3(k)$.
	\end{proof}
	\begin{remark}
		In \cite[Theorem B and Proposition 2.9]{witt-lifting}, the author proves $W_3(k)$-liftability for quasi-$F$-split strict Calabi–Yau varieties in odd characteristic; his construction extends the $W_2(k)$-lifting associated with a chosen quasi-$F$-splitting. Thus, on the common strict projective locus, the existence assertion of Corollary 5.8 is covered by his result. The Hodge--de Rham degeneration used in our $F$-split setting follows from Petrov's theorem \cite[Theorem 1.1]{Pet25}, so no separate degeneration hypothesis is needed.
		
		The emphasis here is the extension criterion of Theorem 5.7: any specified Hodge–Frobenius compatible lifting over $W_n(k)$ extends over $W_{n+1}(k)$, under the stated hypotheses, without assuming projectivity or vanishing of all intermediate structure-sheaf cohomology. The extension need not remain Hodge–Frobenius compatible. This criterion is obtained through the refined period functor developed in the preceding sections, which also applies to the prismatic lifting problem of our Theorem A.
	\end{remark}
	\begin{remark}
		In \cite[Remark 5.7.3]{ST}, the authors raised a question about whether a "canonical lift" (which is just our Hodge–Frobenius compatible lift under certain conditions) always exists. In the following we will give a counterexample that gives a negative answer to the general existence question of such lifts, which also shows that there is a distinction between existing a lift and existing a Hodge–Frobenius compatible lift.
	\end{remark}
	\begin{example}[counterexample of Hodge–Frobenius compatible lift at $W_3\to W_4$]
		Let $S\to\mathbf P^1$ be the rational elliptic surface obtained by resolving the cubic pencil
		\[
		s(x+y)(y+z)(z+x)+txyz=0
		\quad\text{in }\mathbf P^1\times\mathbf P^2.
		\]
		Its singular fibres have types $I_6,I_3,I_2,I_1$, and their parameters remain distinct modulo $7$. Blowing up products of fibre components and the strict transform of the diagonal gives a projective small resolution of $S\times_{\mathbf P^1}S$ with a smooth projective $\mathbf Z_7$-model; take $Y$ to be its base change to $W$. This is the level-six self-product of \cite[Theorem 4.1 and \S5]{SY00}, after the change of base parameter $r=-t/s$ in their equation (10). Let $X=Y\otimes_{W(k)} k$, the logarithmic tangent calculation for its four semistable fibres gives
		\[
		\omega_X\cong O_X,\qquad H^1(X,O_X)=H^2(X,O_X)=0,
		\qquad H^1(X,T_X)=0.
		\]
		Now we show that: $X$ is an $F$-split smooth projective Calabi--Yau threefold $X$ over $k=\overline{\mathbf F}_7$ such that, every $W_m(k)$-lift of $X$ is isomorphic to $Y_m=Y\otimes_W W_m(k)$, and
		\[
		Y_m\text{ is Hodge–Frobenius compatible}\quad\Longleftrightarrow\quad 1\leqslant m\leqslant3.
		\]
		Moreover, all crystalline cohomology groups of $X$ are $W(k)$-free, and all relative Hodge cohomology groups of $Y$ are free, with Hodge--de Rham spectral sequence degenerates.\\
		\\
		\indent First, its Hasse--Witt scalar is $5$ modulo $7$, so $X$ is $F$-split. The cohomological assertions follow from the same Hodge calculation, Deligne--Illusie degeneration, and integral comparison. For the integral input, apply \cite[Theorem 0.3(ii)]{Sch13} to the blow-up at all nodes of the singular self-product over an algebraic closure of $\mathbf Q_7$. This blow-up also dominates the small resolution, whose integral cohomology is a direct summand by the blow-up formula. Then apply \cite[Theorem 0.9]{Min20} for crystalline comparison. In particular, $H^3_{cris}(X/W)$ has rank two. The vanishing of $H^1(X,T_X)$, together with the existence of $Y$, gives uniqueness of the marked lift up to isomorphism at each Witt length.
		\\Now we show the filtered Frobenius obstruction. Put $M=H^3_{cris}(X/W)$ and identify it with $H^3_{dR}(Y/W)$. The slope decomposition and the Hodge numbers give bases and units $\alpha,\beta\in W^\times$ such that
		\[
		M=Wu\oplus Wv,\qquad \varphi(u)=\alpha u,\qquad
		\varphi(v)=7^3\beta v,\qquad
		F^1M=F^2M=F^3M=W(v+cu).
		\]
		Here $\varphi$ is semilinear for Witt Frobenius $\sigma$. Geometric Fontaine--Laffaille theory gives $c\in7^3W$ and the divided maps $\varphi_i=7^{-i}\varphi$ on $F^iM$, compatibly with reduction; see \cite[\S\S4.3--4.4]{LMP23}. The Hodge–Frobenius compatibility condition at length $m$ is
		\[
		\alpha\sigma(c)-7^3\beta c\equiv0\pmod{7^m},
		\]
		which is equivalent to $c\equiv0\pmod{7^m}$.
		\\To see that $c\notin7^4W$, write $Z/\mathbf Q$ for the chosen rational small resolution. Equality of good Euler factors in \cite[Theorem 4.1]{SY00}, together with Chebotarev, identifies $H^3_{\mathrm{\acute et}}(Z_{\overline{\mathbf Q}},\mathbf Q_7)^\vee$, up to global semisimplification, with the global representation $V_f$ attached to
		\[
		f(q)=\eta(q)^2\eta(q^2)^2\eta(q^3)^2\eta(q^6)^2
		=\sum_{r\geqslant1}a_rq^r,\qquad a_5=6,\quad a_7=-16.
		\]
		The reduction of a globally stable lattice is absolutely irreducible: otherwise ordinarity at $7$, squarefree level $6$, and the unipotent inertia at $2$ and $3$ supplied by \cite[Theorem (A)]{Car86} force its semisimplification to be $1\oplus\overline\chi_7^{\,3}$, where $\overline\chi_7$ is the mod-$7$ cyclotomic character. Its trace at $5$ would then be $1+5^3=0$ modulo $7$, contrary to $a_5=6$. If wild inertia at $7$ acted trivially, \cite[Theorem A]{Gee10} would provide a companion form of weight $7+1-4=4$ and level $6$. But the weight-four cusp-form space for $\Gamma_1(6)$ is one-dimensional, also modulo $7$ by tame integral base change, so the normalized companion would equal $f$. Its coefficient at $5$ would have to satisfy $6\equiv5^3\cdot6\equiv1\pmod7$, a contradiction. Absolute irreducibility removes the global semisimplification before restriction to inertia. Thus, for $T=H^3_{\mathrm{\acute et}}(Z_{\overline{\mathbf Q}},\mathbf Z_7)$, the dual $T^\vee$ is a globally stable lattice in $V_f$, and duality transfers nontrivial wild-inertia action to $T/7T$. Put $K=W[1/7]$. Since $K$ is the completion of $\mathbf Q_7^{\mathrm{nr}}$, its Galois group identifies with inertia at $7$, preserving wild inertia. Degree-four torsion-freeness identifies the restricted $T/7T$ with $H^3_{\mathrm{\acute et}}(Y_{\overline K},\mathbf F_7)$.
		\par On the other hand, if $c\in7^4W$, putting $z=v+cu$ gives, for $1\leqslant i\leqslant3$,
		\[
		\varphi_i(z)=7^{3-i}\beta z+
		\bigl(7^{-i}\alpha\sigma(c)-7^{3-i}\beta c\bigr)u.
		\]
		Modulo $7$, these maps split the geometric Fontaine--Laffaille module into two rank-one objects of weights $0$ and $3$. The geometric torsion comparison \cite[Theorem 0.3]{Min20} then splits $H^3_{\mathrm{\acute et}}(Y_{\overline K},\mathbf F_7)$ into two characters. Wild inertia, being pro-$7$, acts trivially on either $\mathbf F_7$-valued character. This contradiction proves $v_7(c)=3$ and hence the asserted failure at $W_4(k)$.
		\par 
	\end{example}
	
\end{document}